\documentclass[leqno,11pt]{article}

\usepackage{amssymb}
\usepackage{euscript}
\usepackage{epic}
\usepackage{eepic}

\evensidemargin\oddsidemargin

\begin{document}

\baselineskip=16pt
\newcounter{index}
\setcounter{index}{1}
\renewcommand{\theequation}{\thesection.\arabic{equation}}  
\newtheorem{theorem}{Theorem}[section] 
\newtheorem{lemma}[theorem]{Lemma}
\newtheorem{proposition}[theorem]{Proposition}
\newtheorem{corollary}[theorem]{Corollary}
\newtheorem{remark}[theorem]{Remark}
\newtheorem{fact}[theorem]{Fact}

\newcommand{\eqnsection}{
\renewcommand{\theequation}{\thesection.\arabic{equation}}
    \makeatletter
    \csname  @addtoreset\endcsname{equation}{section}
    \makeatother}
\eqnsection

\def\r{{\mathbb R}}
\def\e{{\mathbb E}}
\def\p{{\mathbb P}}
\def\z{{\mathbb Z}}
\newcommand{\oo}{\overline}
\def\ee{\mathrm{e}} 
\def\d{\, \mathrm{d}}
\newcommand{\hetzi}{\frac{1}{2}}



\title{Valleys and the maximum local time \\
 for random walk in random environment}

\author{Amir Dembo\thanks{Research partially supported
by NSF grants \#DMS-0406042 and \#DMS-FRG-0244323.}, \
Nina
Gantert\setcounter{footnote}{6}\thanks{Research partially
supported by the DFG and by the European program RDSES}, \ Yuval
Peres\setcounter{footnote}{3}\thanks{Research
supported in part by NSF grants \#DMS-FRG-0244479 and \#DMS-0104073}
\ and Zhan Shi}

\date{} \maketitle


{\leftskip=1.7truecm
\rightskip=1.7truecm
\baselineskip=15pt

\begin{abstract} Let $\xi(n,x)$ be the local time at $x$ for
 a recurrent one-dimensional  random 
 walk in random environment after $n$ steps, and consider the
 maximum $\xi^*(n)=\max_x \xi(n,x)$.
 It is known that $\limsup_n\xi^*(n)/n$ is 
 a positive constant a.s. We prove that 
$\liminf_n(\log\log\log n)\xi^*(n)/n$ is 
a positive constant a.s.; this answers a question of 
P. R\'ev\'esz (1990). The proof is based on an analysis of
the {\em valleys\/} in the environment,
 defined as the potential wells of record depth. 
In particular, we show that
almost surely, at any time $n$ large enough, the random walker has spent 
almost all of its lifetime in the two deepest 
valleys of the environment it has encountered. We also prove a
 uniform exponential tail bound for the
 ratio of the expected total occupation time of a valley and the
 expected local time at its bottom.
\end{abstract}
 \bigskip


 \bigskip


} 

\bigskip
\bigskip

\section{Introduction}
   \label{s:intro}

Let $\omega= (\omega_x)_{x\in \z_+}$ be a collection of i.i.d. random
variables taking values in $(0,1)$. We will denote the distribution of
$\omega$ by $P$. For each $\omega$, we define the random walk in random
environment (RWRE) $(X_n)_{n=0,1,2, \ldots}$ as the Markov chain
taking values in $\z_+$ with 
$X_0=0$ and 
transition probabilities
$P_\omega(X_{n+1} = 1|X_n = 0)=1$, 
$P_\omega(X_{n+1} = x+1|X_n = x) = \omega_x = 1-P_\omega(X_{n+1} =
x-1|X_n = x)$ for $x >0$. For fixed $\omega$, we denote the
distribution of the Markov chain 
$(X_0, X_1, \ldots )$
with $P_\omega$.
As usual, we denote by $\p$ the joint distribution of $(\omega,
(X_n))$. Throughout the paper, we make the following assumptions on the
distribution of the environment $\omega$. Let 
$\rho_i :=
(1-\omega_i)/\omega_i$, 
$i=1,2, \ldots $
\begin{equation}
E (\, \log  \rho_1\,)= \int \, \log \rho_1(\omega) P(d\omega) \, =0\, ,
\label{recurrent}
\end{equation}
\begin{equation}
\hbox{\rm Var}(\, \log  \rho_1\,)
>0 ;
\label{nondet}
\end{equation}
there is $\delta \in (0,1)$ such that 
\begin{equation}
    P(\delta\le \omega_1\le 1-\delta)=1.
    \label{omega0} 
\end{equation}
Assumption (\ref{recurrent}) guarantees that for $P$-almost all
$\omega$, the  Markov chain is recurrent; (\ref{nondet}) excludes the
deterministic case of a simple random walk on the positive integers, and
(\ref{omega0}) is a technical  assumption which could possibly be
relaxed but is used extensively. Usually, one defines in the
same way the RWRE on the integer axis,  but for the questions we will
consider, there is no difference between the two models, so we restrict attention to
the RWRE on the positive integers for simplicity. A key property of recurrent RWRE is
its strong localization: under our assumptions, Sinai \cite{sinai} showed that
that $X_n/(\log n)^2$  converges in distribution.
A lot more is known about this model; we refer to the  survey 
by Zeitouni~\cite{zeitouni} for 
limit theorems, large deviations results, and for further references. 

 Let $\xi(n,x): = |\{0 \le j \le n: X_j =x\}|$
denote the local time of the RWRE in $x$ at time $n$ and  $\xi^*(n) :=
\sup_{x\in \z_+} \xi(n,x)$ the maximal local time at time $n$. 
It was proved in \cite{mama} that for each non-decreasing function $\varphi$, $\limsup_{n\to \infty} \, {\xi^*(n)\over \varphi(n)}$ and $\liminf_{n\to \infty} \, {\xi^*(n)\over \varphi(n)}$ are $\p$-almost surely (possibly degenerate) constants. For the $\limsup$ behavior of $\xi^*(n)$, it 
was
shown  in \cite{r90} and \cite{zhan1} that
$$ 
\limsup_{n\to \infty} \, {\xi^*(n)\over n} >0
    \qquad \hbox{\rm $\p$-a.s.}
$$
(Clearly this $\limsup$ is at most $1/2$.)
In his book, R\'ev\'esz~\cite{r90} raised the problem of determining the 
$\liminf$ behavior of $\xi^*(n)$.
\noindent Our main result is the following.

\medskip

\begin{theorem}
 \label{t:main} 
 There exists a constant $c\in (0,\infty)$ such that
 \begin{equation}
    \liminf_{n\to \infty} \, {\xi^*(n)\over n/\log \log\log n} =
    c, \qquad \hbox{\rm $\p$-a.s.}
    \label{main}
 \end{equation}
\end{theorem}

\smallskip

In particular, (\ref{main}) disproves the conjecture on page 303 of
R\'ev\'esz \cite{r90}. We will shortly give a heuristic argument
which explains why the three logarithms appear.

The {\em potential\/} corresponding to the RWRE is $V(x)=\sum_{i=1}^x \log \rho_i$,
${x\in \z_+}$.
As is well known, the potential governs the behavior of the RWRE in several senses, e.g. 
\begin{itemize}
\item 
In an excursion starting from any site $b$, the logarithm of 
the expected number of visits to a site $x$ before returning to $b$
is roughly the potential difference $V(b)-V(x)$, see (\ref{expu}).
\item Starting from the origin,
the logarithm of the expected time to reach a site $x$ is roughly $\max_{y \le x} V(y)$;
 see (\ref{l:hitting}) for an upper bound, and observe that a
 similar lower bound follows from  (\ref{expu}). 
\end{itemize}

The proof of Theorem \ref{t:main} is based on an
analysis of the {\em valleys} in the potential, which is of independent interest.
By a ``valley'' we mean a potential well of record depth; 
see \S \ref{subs:valley} for a precise definition.

We will partition the environment into valleys, and show that at any
time $n$, the particle performing RWRE has almost surely spent almost all
of its lifetime in the two deepest valleys it has encountered.
This almost
sure localization theorem (Theorem \ref{t:spend}) can be considered as
the second main result of the paper. 
Furthermore, we define in (\ref{Lambda}) the {\em effective width} of a valley as the ratio of the expected total occupation time of the valley and the 
expected local time at its bottom, and prove a uniform exponential tail bound 
(\ref{expo}) for the effective width of valleys.
The reason for the term  ``effective width'' is that 
most of the occupation time in a valley is spent at sites
where the potential is within an additive constant from its minimum in the valley;
the number of these sites is the effective width, up to a multiplicative constant.

Theorem \ref{t:main} is then established  as follows:
    
Due to scaling properties of the potential, the 
depths of 
successive 
valleys grow at a geometric rate, 
whence the distance between bottoms of successive valleys also 
exhibit geometric growth,
resulting with $O(\log R)$ valleys in a large interval $[0,R]$. By time 
$n$ the random walker reaches a distance of order  
$(\log n)^2$ from the origin, thus
visiting an order of $\log\log n$ valleys. 
The exponential tail bounds on effective widths imply that a.s., 
for all $k$, the
$k$th valley encountered has effective width at most $O(\log k)$;
conversely, a.s. for infinitely many $k$ the  effective width is at least $c\log k$.
Hence, a.s. the maximal effective width of  
valleys seen by the walker up to time $n$ is at most of order 
$\log \log\log n$, and up to a constant factor,
this effective width is realized infinitely often.

The 
paper is organized as follows. In Section
\ref{s:valley}, we introduce the notion of valleys, and describe some
scaling properties of such valleys. Section \ref{s:particle} is devoted
to the study of the behavior of the RWRE within the valleys. We first
give some background on hitting times and excursions. We then compare
the occupation time of different valleys and prove that the RWRE spends
most of its time in the last two visited valleys: Theorem \ref{t:spend}
is the main result of this section. In Section \ref{s:lt}, we compare
the occupation time of valleys with the local time in sites. Our main
tool here is to average over excursions of the RWRE. 
This comparison motivates our definition of the \lq\lq
effective 
width\rq\rq\ of the valleys, whose asymptotic growth is studied in  
the second part of Section \ref{s:lt}. Similarly to 
Section \ref{s:valley} this part does not concern 
the random walk, but only the environment. 
Finally, Theorem \ref{t:main} is proved in Section \ref{s:proof}.

\section{Valleys}
\label{s:valley}

Recall that the potential
$V$ is a function of the environment, defined as follows:
$$
V(x) :=
\cases{ 
\sum_{i=1}^x \log \rho_i, &$x=1,2, \ldots $, \cr\cr
        0, &$x=0$. \cr}
$$
Note that $V$ is itself a sum of i.i.d. random variables, which are
bounded by 
$C:= |\log \delta - \log (1-\delta)|$, see (\ref{omega0}). For fixed $\omega$, $P_\omega$ is a reversible Markov chain, hence an electrical network in the sense of \cite{doylesnell}.
The conductance of the bonds is
\begin{equation} 
C_{(x, x+1)} = \ee^{-V(x)}, \, x = 0,1,2, \ldots
    \label{conductance} 
\end{equation}
and the reversible measure $\mu$ 
(which is unique up to multiplication by a constant), is given by
\begin{equation} 
    \mu(x) =
    \cases{ \ee^{-V(x)} + \ee^{-V(x-1)}, &$x=1,2, \ldots $, \cr\cr
        1, &$x=0$. \cr}\, 
    \label{mu}
\end{equation}
For background on reversible Markov chains, we refer to \cite{doylesnell}.

\subsection{Definition of valleys}
\label{subs:valley}
Fix a constant $K_0>0$. We set $\theta_0:=0$ and 
\begin{eqnarray*}
    \eta_0 
 &:=& \inf \Big\{ i > 0: \; V(i) - \min_{0\le j\le i} V(j) \ge
    K_0 \Big\} ,
    \\
    b_0
 &:=& \sup \Big\{ i < \eta_0: \; V(i) = \min_{0\le j\le \eta_0} V(j)
    \Big\}.
\end{eqnarray*}
We now define, for $k \ge 1$, inductively:
\begin{eqnarray*}
    \theta_k
 &:=& \inf \Big\{ i> \eta_{k-1}: \; V(i) \le V(b_{k-1}) \Big\}, 
    \\
    H_{k-1}^+
 &:=& \max_{\eta_{k-1} \le j\le \theta_k} V(j) - V(b_{k-1}),
    \\
    \eta_k
 &:=& \inf \Big\{ i> \theta_k: \; V(i) - \min_{0\le j\le i} V(j)
    \ge H_{k-1}^+ \Big\},
    \\
    b_k
 &:=& \sup\Big\{ i<\eta_k: \; V(i) = \min_{\theta_k \le j\le
    \eta_k} V(j) \Big\},
    \\
    H_k^-
 &:=& \max_{\eta_{k-1} \le j\le \theta_k} V(j) - V(b_k).
\end{eqnarray*}
Let now
\begin{equation}
    m_k:= \inf \, \Big\{ i> \eta_{k-1}: \; V(i) =
    \max_{ \eta_{k-1} \le j\le \theta_k} V(j) \Big\} 
    \,.
    \label{mk}
\end{equation}  

The piece $(V(i), \; m_k\le i < m_{k+1})$ is the $k$-th valley, $H_k^-$
the left height of this valley, and $H_k^+$ the right height. We call
$$
H_k := \min\left\{ H_k^-, \; H_k^+ \right\},
$$

\noindent the height of the $k$-th valley. Also, $b_k$ is called
the bottom of the $k$-th valley. 

\bigskip

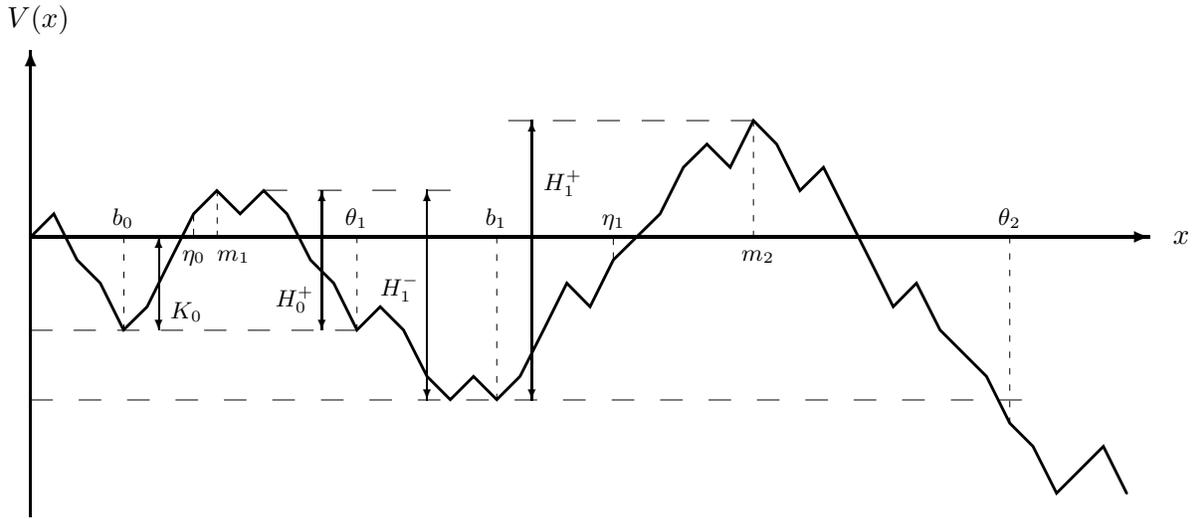
\begin{center}
\unitlength = 3.1mm 
 \begin{picture}(53,25)(-1,-14)
\thicklines 
\put(0,0){\vector(1,0){48}}  
\put(0,-12){\vector(0,1){20}}
\put(49,-0.3){$x$}\put(-1,9){$V(x)$}
\Thicklines 
\path(0,0)(1,1)(2,-1)(3,-2)(4,-4)(5,-3)(6,-1)(7,1)(8,2)(9,1)(10,2)(11,1)(12,-1)
(13,-2)(14,-4)(15,-3)(16,-4)(17,-6)(18,-7)
(19,-6)(20,-7)
(21,-6)(22,-4)(23,-2)(24,-3)(25,-1)(26,0)(27,1)(28,3)(29,4)(30,3)
(31,5)(32,4)(33,2)(34,3)(35,1)(36,-1)(37,-3)
(38,-2)(39,-4)(40,-5)(41,-6)(42,-8)(43,-9)(44,-11)(45,-10)(46,-9)(47,-11)

\thinlines
\dashline{0.3}(4,-4)(4,0)\put(3.5,0.5){\footnotesize$b_0$}
\dashline{0.8}(7,1)(7,0)\put(6.5,-1){\footnotesize$\eta_0$}
\dashline{0.3}(8,2)(8,0)\put(8,-1){\footnotesize$m_1$}
\dashline{0.3}(14,-4)(14,0)\put(13.5,0.5){\footnotesize$\theta_1$}
\dashline{0.3}
(20,-7)(20,0)\put(19.5,0.5){\footnotesize$b_1$}
\dashline{0.8}(25,-1)(25,0)\put(24.5,0.5){\footnotesize$\eta_1$}
\dashline{0.3}(31,5)(31,0)\put(30.5,-1){\footnotesize$m_2$}
\dashline{0.3}(42,-8)(42,0)\put(41.5,0.5){\footnotesize$\theta_2$}
\put(5.5,0){\vector(0,-1){4}}\put(5.5,-4){\vector(0,1){4}}\put(6,-3.5){\footnotesize$K_0$}
\put(12.5,2){\vector(0,-1){6}}\put(12.5,-4){\vector(0,1){6}}\put(10.5,-3){\footnotesize$H_0^+$}
\put(17,2){\vector(0,-1){9}}\put(17,-7){\vector(0,1){9}}\put(15,-2.5){\footnotesize$H_1^-$}
\put(21.5,5){\vector(0,-1){12}}\put(21.5,-7){\vector(0,1){12}}\put(22,2){\footnotesize$H_1^+$}
\dashline{1}(0,-4)(14,-4)
\dashline{1}(0,-7)(42.5,-7)
\dashline{1}(10,2)(18,2)
\dashline{1}(20.5,5)(31,5)
\end{picture}\end{center}

\bigskip

\noindent{\it Remark.}
(i) In words, $m_k$ is the beginning of the $k$-th valley. Note that $(\theta_k)_{k\ge 0}$ and $(\eta_k)_{k\ge 0}$ are
sequences of stopping times (with respect to the natural
filtration of the potential $V$), whereas $(b_k)_{k\ge 1}$ and $(m_k)_{k\ge 1}$ are not.

(ii) Our definition of valleys is not exactly the standard
definition of valleys in the sense of Sinai \cite{sinai}.
However, it follows from our definition that almost surely the heights
$(H_k, \; k\ge 1)$ are increasing.

(iii)
Here is a (very) rough description of the asymptotic behavior of the RWRE. When $k$ is large,
the time needed for the RWRE to exit from the $k$-th valley is of order $\ee^{H_k}$ (see (\ref{l:hitting}) and Lemma \ref{l:exit})); and since $H_k$ is of order $\ee^k$ (Lemma \ref{l:Hk}), we have: $n \approx \ee^{H_{N_n}}$, where $N_n$ is the number of valleys visited by the RWRE in the first $n$ steps. This leads to: $H_{N_n} \approx \log n$. On the other hand, $V$ being the partial sum process of i.i.d.\ bounded mean-zero random variables, $H_k \approx x_k^{1/2}$ for any site $x_k$ in the $k$-th valley. Therefore, $x_{N_n}$ is of order $(\log n)^2$; i.e., the maximal distance to the origin of the RWRE in the first $n$ steps is of order $(\log n)^2$. In fact, a famous result of Sinai~\cite{sinai} says that ${X_n\over (\log n)^2}$ converges in distribution (under $\p$) to a non-degenerate limit.

\subsection{Heights of valleys}
\label{subs:height}

We now consider the asymptotic growth of the heights of the valleys.

\medskip

\begin{lemma}
 \label{l:Hk} 
 We have, $P$-almost surely,
 \begin{equation}
     \log H_k \sim \log H_k^+ \sim \log H_k^- \sim k, \qquad k\to
     \infty.
     \label{logHk1}
 \end{equation}
\end{lemma}

\medskip

\noindent {\it Proof.} Assume for a moment that $V$ is a Brownian motion. 
Then, the strong Markov property
at $\theta_k$
and scaling properties imply that
$({H_{k-1}^+\over H_k^+}, \; k\ge 2)$ is a sequence of i.i.d. random
variables  with common uniform distribution 
on 
$(0,1)$. In particular,
$E (\log {H_k^+\over H_{k-1}^+} )=1$. More precisely, since
$(\eta_k, k\geq 1)$ is a sequence of stopping times, the random
variables 
${H_{k-1}^+\over H_k^+}$, $k\ge 2$,
are
independent,
 and the probability of the event 
$\{H_k^+ \geq (1+c)H_{k-1}^+\}$ is the probability that a standard
Brownian motion hits $c$ before hitting $-1$. By the
law of large numbers, $P$-almost surely,
\begin{equation}
    \log H_k^+ = \log H_1^+ + \sum_{i=2}^k \log {H_i^+ \over
    H_{i-1}^+} \sim k, \qquad k\to \infty.
    \label{logHk}
\end{equation}
\noindent Further, the strong Markov property at the stopping time $\theta_k$
implies that ${H_k^- - H_{k-1}^+ \over
H_{k-1}^+}$ is an exponential random variable with mean 1.
More precisely, the conditional distribution of ${H_k^- - H_{k-1}^+ \over
H_{k-1}^+}$, given $H_{k-1}^+ = a$, is the distribution of
$| \inf_{t < \sigma (a)} B_t |\cdot a^{-1}$, where $(B_t)$ is a standard
Brownian motion and
$\sigma(a): = \inf\{s: B_s -\inf_{u < s} B_u  =a\}$. By scaling, this
distribution does not depend on $a$, hence equals the distribution of
$| \inf_{t < \sigma(1)} B_t |$. L\'evy's identity tells us that $(B_t
-\inf_{s < t} B_s, |\inf_{s < t} B_s|)$ has the same distribution as
$(|B_t|, L_t)$ where $(L_t)$ is the local time of $(B_t)$ at $0$
(c.f. \cite[Theorem VI.2.3]{RY}).
Therefore, $| \inf_{t < \sigma (1)} B_t |$ has the same distribution as
$L_\tau$, where $\tau:= \inf\{t: |B_t| =1\}$, and the 
distribution of $L_\tau$ is known to be exponential with mean $1$
(for example, see Formula 3.3.2, page 213 of \cite{borodin}).
Using the Borel--Cantelli lemma, we see that $P$-almost surely,
\begin{equation}
    \log \Big( {H_k^- \over H_{k-1}^+} \Big)  = O(\log\log k), 
    \qquad k\to \infty,
    \label{Hk/Hk-1}
\end{equation}

\noindent and thus (\ref{logHk}) yields $\log H_k^- \sim k$,
$P$-almost surely. This would prove the lemma if $V$ 
was 
a Brownian motion.

In our case, $V$ is the partial sum process associated with a
sequence of i.i.d. bounded mean-zero random variables, so we have to be
more careful. Let $k\ge 1$. We look at the random walk $(V(i+
\theta_k)-V(\theta_k), \; i\ge 0)$, which is independent of $(V(i),
\; i\le \theta_k)$ (thus of $H_{k-1}^-$ and $H_{k-1}^+$). This
random walk can be embedded into a Brownian motion, say $(B_k (t), \;
t\ge 0)$, in the sense of Skorokhod embedding,
making $V(i+\theta_k)- V(\theta_k)=B_k (t_i), \; i\ge 0$, a random sequence of
points on the path of $t \mapsto B_k (t)$, such that the maximum  
of the height differences $|B_k (t)-B_k (t_i)|$ for $t \in [t_i,t_{i+1}]$ is
at most $C$. For any $r>0$, let
$$
\sigma_k (r):= \inf\{ t>0: \, B_k (t) - \inf_{s\in [0,t]} B_k (s) = r\}\,,
$$
and 
$$
{\widetilde H}_k^- (r) := r + |\inf_{0\le t\le \sigma_k (r)}\, B_k (t)|\,.
$$ 
Note that with $V(b_{k-1}) - C \leq V(\theta_k) \leq V(b_{k-1})$,
given $H_{k-1}^+ = a>0$, we have that
$$
a+V(\theta_k)-V(b_k) \leq H_k^- \leq a + C + V(\theta_k)-V(b_k) \,.
$$
Further, $V(b_k)-V(\theta_k)$ is the minimum of $B_k (t_i)$ for
those $i$ such that $t_i \in [0,t_{\eta_k-\theta_k}]$, and 
since the Brownian increments between $B_k (t_i)$ are 
of height at most $C$, we have that
$$
V(b_k)-V(\theta_k) - C 
\leq \inf_{0 \leq t \leq t_{\eta_k-\theta_k}} B_k(t)
\leq V(b_k)-V(\theta_k) \,.
$$
We thus conclude that if $H_{k-1}^+ =a>2C$, then
\begin{equation}
    {\widetilde H}_{k}^- (a - 2C)  
\le H_k^- \le {\widetilde H}_{k}^- (a+2C)\,.
    \label{HkHkHk}
\end{equation}
More precisely, 
by the time $\sigma_k (a+2C)$ the Brownian 
motion made an increment of 
$a+2C$ over its minimal 
value and by the time $\sigma_k (a-2C)$ it made
an increment of 
$a-2C$ over its minimal value. Since
$\eta_k-\theta_k$ corresponds to the first value of $i$
where $B_k(t_i)$ makes an increment of at least $a$ from 
its minimum, and the Brownian increments between the points
$B_k(t_i)$ are at most of height $C$, a fortiori, 
$$
\sigma_k (a-2C) \leq t_{\eta_k-\theta_k} \leq 
\sigma_k (a+2C)\,, 
$$
which by the monotonicity of $u \mapsto \inf_{0 \leq t \leq u} B_k (t)$
yields the inequality (\ref{HkHkHk}). 

Similarly, we embed the
random walk $(V(j+ \eta_k)-V(\eta_k), \; j\ge 0)$ as a random
sequence of points $W_k(s_j)$ on the path of an independent Brownian motion
denoted $(W_k(s), \; s\ge 0)$, 
such that the maximum of the height differences
$|W_k(s)-W_k(s_j)|$ for $s \in [s_j,s_{j+1}]$ is at most $C$,
and without loss of generality,
we assume that we are still working on the same probability space. 
Note that $V(\eta_k)-V(b_k)$ is within 
distance
$C$ of $H_{k-1}^+$ and that
$$
H_k^+ = \max_{0 \leq j \leq \theta_{k+1}-\eta_k} W_k(s_j) + V(\eta_k)-V(b_k) \,,
$$
where $\theta_{k+1}-\eta_k$ corresponds to the first value of $j$
such that $W_k(s_j) \leq V(b_k)- V(\eta_k)$. Therefore, 
by a similar line of reasoning as before,
given $H_{k-1}^+ = a> 2C$, we have that 
$$
S_k (-(a-2C)) \leq s_{\theta_{k+1}-\eta_k} \leq S_k (-(a+2C))\,,
$$
where $S_k (r):= \inf\{ s \geq 0: \, W_k (s) = r\}$.
Consequently, then also
\begin{equation}
    {\widetilde H}_k^+ (a-2C)
\le H_k^+ \le {\widetilde H}_{k}^+ (a+2C) \,,
    \label{HkHkHk1}
\end{equation}
where for any $r>0$,
$$
{\widetilde H}_{k}^+ (r) := r
+ \sup_{0\le s\le S_k (-r)}W_k (s) \,.
$$

Recall that $H_k^+ = V(m_{k+1})-V(b_k)$, is non-decreasing, and further 
$$
H_k^+ - H_{k-1}^+ \geq V(b_{k-1}) - V(b_k) \geq V(\theta_k)-V(b_k) \,,
$$
which is non-negative, and dominates the law of the negative part of
$\log \rho_0$. Thus, by (\ref{nondet}) we see that
$H_k^+\to \infty$, $P$-almost surely. Fixing 
$\varepsilon > 0$, we thus have that 
$P$-almost surely, $\varepsilon H_{k-1}^+ \geq 2 C$ for all $k$ large enough,
in which case we have from (\ref{HkHkHk}) and (\ref{HkHkHk1}) that  
\begin{equation}\label{eq:amir-bd}
{\widetilde H}_{k}^{\pm} ( (1-\varepsilon) H_{k-1}^+)  
\le H_k^{\pm} \le {\widetilde H}_{k}^{\pm} ((1+\varepsilon) H_{k-1}^+ )\,.
\end{equation}
Without loss of generality we take the Brownian motions
$B_k(\cdot)$, $W_k(\cdot)$, $k=1,2,\ldots$, to be independent, 
and consequently, so are ${\widetilde H}_{k}^{\pm} (\cdot)$. Further, 
by the scaling properties of the Brownian motion, the law of
$r^{-1} {\widetilde H}_{k}^{\pm} (r)$ is independent
of $r>0$ and $k$, resulting with i.i.d. random variables  
$$
Z_k^{\pm} := {{\widetilde H}_{k}^{\pm} ( u H_{k-1}^+) \over u H_{k-1}^+}\,,
$$   
whose law
is independent of $u>0$. As we have already seen,
$-1 + Z_k^-$ has the exponential distribution of mean $1$ (being
the same as $|\inf _{t<\sigma(1)} B_t|$) while $1/Z_k^+$ has the uniform law
on $(0,1)$. Consequently, $E(\log 
Z_1^+) = 1$ and
$$
k^{-1} \sum_{i=1}^k \log Z_i^+ \to 1 
$$
$P$-almost surely. Since (\ref{eq:amir-bd}) holds for all but finitely 
many values of $k$ and $\log (1 \pm \varepsilon)$ can be arbitrarily small, 
it follows that also
$$
k^{-1} \sum_{i=2}^k \log {H_i^+ \over H_{i-1}^+} \to 1\,,
$$
$P$-almost surely. That is, $\log H_k^+ \sim k$, $P$-almost surely.

A Borel--Cantelli argument as in the proof of
(\ref{Hk/Hk-1}), using (\ref{eq:amir-bd}), easily implies that $\log (
{H_k^-/H_{k-1}^+})  = O(\log\log k)$, $P$-almost surely. Thus $\log
H_k^- \sim k$, $P$-almost surely. This completes the proof of
Lemma \ref{l:Hk}.\hfill$\Box$

\medskip

\begin{lemma}
 \label{l:Hk-Hk-1}
 Let $\varepsilon>0$. We have, $P$-almost surely for all
 sufficiently large $k$,
 \begin{equation}
     H_k-H_{k-1}^+ \ge (H_{k-1}^+)^{1-\varepsilon} .
     \label{Hk-Hk-1}
 \end{equation}
\end{lemma}

\medskip

\noindent {\it Proof.} Observe that
$$
P\Big( {H_k\over H_{k-1}^+} < 1+ \ee^{-\varepsilon k/2}\Big) \le
P\Big( {H_k^+\over H_{k-1}^+} < 1+ \ee^{-\varepsilon k/2}\Big) +
P\Big( {H_k^-\over H_{k-1}^+} < 1+ \ee^{-\varepsilon k/2}\Big).
$$

\noindent The distributions of ${H_k^+\over H_{k-1}^+}$  and
${H_k^-\over H_{k-1}^+}$ have already been mentioned in
the case of a Brownian potential $V$: ${H_{k-1}^+\over H_k^+}$ is uniformly
distributed 
on 
$(0,1)$, whereas ${H_k^- - H_{k-1}^+ \over
H_{k-1}^+}$ is an exponential random variable with mean 1.
Therefore, $\sum_k P\big( {H_k^+\over H_{k-1}^+} < 1+
\ee^{-\varepsilon k/2} \big) <\infty$ and $\sum_k P\big( {H_k^-\over
H_{k-1}^+} < 1+ \ee^{-\varepsilon k/2} \big) <\infty$. As a consequence,
$\sum_k P\big( {H_k\over H_{k-1}^+} < 1+ \ee^{-\varepsilon k/2} \big)
<\infty$. 

For our partial sum potential, we can easily use 
(\ref{eq:amir-bd})
to see that $\sum_k P\big( {H_k\over H_{k-1}^+} < 1+
\ee^{-\varepsilon k/2}\big) <\infty$ still holds. By the
Borel--Cantelli lemma, $P$-almost surely for $k$ large enough, $H_k -
H_{k-1}^+ \ge H_{k-1}^+\ee^{-\varepsilon k/2}$. This yields
(\ref{Hk-Hk-1}), as we know from Lemma \ref{l:Hk} that $\log H_{k-1}^+
\sim k$, $P$-almost surely.\hfill$\Box$

\subsection{Other facts about valleys}
\label{subs:otherfacts}

Throughout the paper, we will subsequently use some 
asymptotic properties 
of the valleys. 
First, note that
$$
K_0 + \sum_{i=1}^k (H_i^- + H_i^+) \ge \max_{0\le x, \, y \le
m_{k+1}} |V(x)-V(y)| \ge 
\max_{x\in [0,\, m_{k+1}]} |V(x)| \ge 
\frac{1}{2} H_k^-.
$$
Hence, with $m_k \to \infty$, 
applying Chung's law of the iterated
logarithm for the potential $V$, we   
have
for each $\varepsilon \in (0,1/4)$, that $P$-almost surely for all
sufficiently large $k$,
\begin{equation}
\label{mk1>Hk2}
K_0 + \sum_{i=1}^k (H_i^- + H_i^+) \ge
m_{k+1}^{(1-0.5 \varepsilon)/2} \ge (H_k^-)^{1-\varepsilon}.
\end{equation}
In view of Lemma \ref{l:Hk} 
the first inequality in (\ref{mk1>Hk2})
implies that $P$-almost surely,   
\begin{equation}
    m_k \le b_k \le 
m_{k+1} \le 
H_k^{2+\varepsilon},
    \label{bk<Hk2}
\end{equation}
for all sufficiently large $k$. Further, by the same reasoning 
we have that $P$-almost surely,
\begin{equation}
\label{loglogm_k}
\log\log m_k \sim \log k, \qquad \hbox{ \rm for } k \to \infty\, .
\end{equation}

We will also make use of the following: for each $\varepsilon \in
(0,1)$,  we have $P$-almost surely for all $k$ large enough,
\begin{eqnarray}
    \max_{m_k\le y \le z < b_k} (V(z)-V(y))
 &\le& H_{k-1}^+ - (H_{k-1}^+)^{1-\varepsilon},
    \label{depth1}
    \\
    \max_{b_k\le y \le z < m_{k+1}} (V(y)-V(z))
 &\le& H_k^+ - (H_k^+)^{1-\varepsilon}.
    \label{depth2}
\end{eqnarray}

\noindent Moreover, 
\begin{eqnarray}
    \min_{x\in [\eta_k, \, m_{k+1})} V(x) 
 &\ge& V(b_k)+ (H_{k-1}^+)^{1-\varepsilon} ,
    \label{depth3}
    \\
    \max_{b_k\le y \le z < \eta_k} (V(y)-V(z))
 &\le& H_{k-1}^+ - (H_{k-1}^+)^{1-\varepsilon}.
    \label{depth4}
\end{eqnarray}

We next outline the proof of (\ref{depth1}) in case $V$ is
a Brownian motion. A similar argument as in the
proof of Lemma \ref{l:Hk} will then confirm that
(\ref{depth1}) holds also when $V$ is a partial sum process. 
With $H_{k-2}^+$ measurable on the stopped $\sigma$-field at 
$\theta_{k-1}<\eta_{k-1}$, for
$V(\cdot)$ a Brownian motion we have by the strong Markov property
at $\eta_{k-1}$ that conditionally on $H_{k-2}^+ = a>0$
the process $U(s):= (V(s+\eta_{k-1}) - V(\eta_{k-1})+a, \,
0\le s \le \theta_k - \eta_{k-1})$ is also a Brownian motion,
starting from $U(0)=a$ and killed upon first hitting $0$
(at time 
$\theta_k-\eta_{k-1}=:S(0)$). 
Of course, in this case also
$H_{k-1}^+ = \sup_{0\le s \le S(0)} U(s) =: H$
and
$m_k-\eta_{k-1} = \inf \{ s \geq 0 : U(s) = H \} =: m_H$.
Thus, denoting by $P_x$ the probability law of a Brownian motion $U(\cdot)$
starting at $U(0)=x$ and by $S(y) := \inf \{ t \geq 0 : U(t) = y \}$
the corresponding first hitting time of $y$, it follows that for
any $a > 0$ and $k \geq 2$,
$$
P \Big( \max_{m_k \le y \le z<\theta_k} (V(z)-
    V(y)) > H_{k-1}^+ - (H_{k-1}^+)^{1- \varepsilon} 
    \, | \, H_{k-2}^+ = a \Big)
\leq P_a (H < \lceil a \rceil) + \sum_{h=\lceil a \rceil}^\infty J(a,h)
$$
where for integer $h \geq 1$,
$$
J(a,h) =
P_a \big( h \leq  H < h+1, \max_{m_H \le y \le z< S(0)} (U(z)- U(y))
> H - H^{1- \varepsilon} \big) \,.
$$
Since $H \geq U(z)$ and $U(y) \geq 0$, the
event whose probability is $J(a,h)$ requires the existence of random
times $m_H<y<z<S(0)$ with $U(m_H) \geq h$,
$U(y)<(h+1)^{1-\varepsilon} =: u$,
$U(z)> h - (h+1)^{1-\varepsilon} =: v$
and $U(S(0)) = 0$, while
$0<U(s)<h+1$ for all $s < S(0)$. It is easy to see that
$h > 2 (h+1)^{1-\varepsilon}$ for any $h \geq a \geq 3^{1/\varepsilon}$,
in which case by continuity of the Brownian path and the preceding reasoning,
\begin{eqnarray*}
J(a,h)&\leq&
P_a (S(h)<S(0)) P_h (S(u) < S(h+1)) P_{u} (S(v) < S(0))
P_{v} (S(0)<S(h+1)) \\
&=& \frac{a u (h+1-v)}{h (h+1-u) v (h+1)}
\leq  8 a  h^{-(2+2\varepsilon)} \,.
\end{eqnarray*}
Hence, $\sum_{h \geq a} J(a,h) \leq C a^{-2\varepsilon}$ for
a finite constant $C=C(\varepsilon) \geq 1$. Further,
$P_a (H < \lceil a \rceil) \leq a^{-1}$, so we conclude that
$$
P \Big( \max_{m_k \le y \le z<\theta_k} (V(z)- V(y))>
H_{k-1}^+ - (H_{k-1}^+)^{1- \varepsilon} \Big) \le
P (H_{k-2}^+ \leq 3^{1/\varepsilon}) + 2 C E ( (H_{k-2}^+)^{-\varepsilon} )
$$
which is summable in $k$ (recall that $H_1^+ \geq K_0$ and
$H_{i-1}^+/H_i^+$, $i \geq 2$, are i.i.d. uniform $(0,1)$ random variables).
Thus, $P$-almost surely for all large $k$,
$$
\max_{m_k \le y \le z<\theta_k} (V(z)- V(y))
\le H_{k-1}^+ - (H_{k-1}^+)^{1- \varepsilon} \,.
$$
A similar (and easier) argument shows that,
$P$-almost surely for all large $k$
$$
\max_{\theta_k \le y \le z<b_k} (V(z)-V(y)) \le
H_{k-1}^+ - (H_{k-1}^+)^{1-\varepsilon} \,,
$$
yielding (\ref{depth1}) when $V$ is a Brownian motion.

The proof of (\ref{depth2}) is very similar.
The proofs of (\ref{depth3}) and (\ref{depth4}) are even easier
since $H_{k-1}^+$ is measurable on the stopped $\sigma$-field at
$\eta_k$ and $\theta_k$, which is where
we apply the strong Markov property when proving
(\ref{depth3}) and (\ref{depth4}), respectively.

\section{Particle in the valleys}
\label{s:particle}

In this section, we will consider the RWRE and give estimates on
hitting times, exit times and excursions.

\subsection{Hitting time}
\label{subs:hitting}

For any $x\in \z_+$, define
$$
T(x) := \inf\left\{ n
\geq 1
: \; X_n = x \right\},
$$
the first hitting time of $x$ by the particle. 
The inequality \cite[(A.1)]{golosov} states that 
for any 
$x \geq 1$,
\begin{equation}
 \label{l:hitting} 
    E_\omega (\, T(x)\,) \le x^2 \, \exp\Big( \, \max_{0\le i\le
    j< x} (V(j)-V(i)) \Big) .
\end{equation}
A consequence of (\ref{l:hitting}) is that for any $k\ge 2$ and any
$\lambda \ge 1$,
\begin{equation}
    P_\omega\left( T(b_k) \ge \lambda\right) \le {b_k^2\over \lambda}
    \, \ee^{H_{k-1}^+}.
    \label{cheby}
\end{equation}

Another result we will be frequently using concerns the almost sure
asymptotic behavior of $T(x)$ when $x\to \infty$. The following is a
consequence of the law of the iterated logarithm for RWRE, stated in
Theorems 27.8 and 27.9 of R\'ev\'esz \cite{r90}.

\medskip

\begin{fact}
 \label{f:revesz} {\bf (R\'ev\'esz \cite{r90})} We have,
$$
    \lim_{x\to \infty} \, {\log\log T(x) \over \log x} = {1\over 2} ,
    \qquad \hbox{\rm $\p$-a.s.}
$$ 
\end{fact}

\medskip

Consider the $k$-th valley $(V(i), \; m_k\le i
< m_{k+1})$. Let a
particle $(X_n, \; n\ge 0)$ start from the bottom $X_0=b_k$ of the
valley. We are interested in
$$
\tau_k := \inf\left\{ n>0: \; X_n \notin (m_k, m_{k+1})\right\},
$$

\noindent the first exit time of the particle from the valley.

\medskip

\begin{lemma}
 \label{l:exit} 
For some $c_0<\infty$, 
any $k\ge 1$ and $m\ge 1$,
\begin{equation}
  P_\omega\left( 
\tau_k < m \, | \, X_0= b_k
\right) \le 
c_0 m \ee^{-H_k}\,.
    \label{exit}
 \end{equation}
\end{lemma}

\medskip

\noindent {\it Proof.} Considering the side from which the particle
exits the valley, we see that
$$
P_\omega( \tau_k < m
\, | \, X_0= b_k
) \le P_\omega(T(m_k)<m
\, | \, 
X_0=b_k) +
P_\omega(T(m_{k+1})<m 
\, | \, 
X_0=b_k) \,,
$$
hence (\ref{exit}) is 
just
a consequence of \cite[Lemma 7]{golosov},
the definition of $H_k$,
and the fact that increments of $V$ are bounded by $C$.
\hfill$\Box$

\medskip

\begin{corollary}
 \label{c:exit} 
 For any $k\ge 1$ and $a>0$,
 $$
    E_\omega\left( \ee^{-a \tau_k} 
\, | \, X_0=b_k
\right) \le {2 c_0\ee^{-a} \over
    (1- \ee^{-a})} \, \ee^{- H_k}\leq \frac{2 c_0}{a}\ee^{- H_k}\, .
 $$
\end{corollary}

\medskip

\noindent {\it Proof.} By changing the order of summation,
$$
E_\omega\left( \ee^{-a \tau_k} 
\, | \, X_0=b_k
\right) =
(1-\ee^{-a})\sum_{m=1}^\infty \ee^{-a m} P_\omega\left( \tau_k \leq m 
\, | \, X_0=b_k
\right) .
$$
Replacing $P_\omega\left( \tau_k \leq m 
\, | \, X_0=b_k
\right)$ by
$P_\omega\left( \tau_k < m+1 
\, | \, X_0=b_k
\right)$ and using (\ref{exit}), the
corollary follows easily. \hfill$\Box$

\bigskip

We note for further reference that for $b < x < i$,
\begin{equation}
\label{hitbebeforei}
    P_\omega\left(T(b) < T(i) | X_0 = x\right) = 
    \sum\limits_{j=x}^{i-1}\ee^{V(j)}
    \Big(\sum\limits_{j=b}^{i-1}\ee^{V(j)}\Big)^{-1} .
\end{equation}
This follows from direct computation, using (\ref{conductance}), 
see also \cite{zeitouni}, formula (2.1.4).

\subsection{Excursions}
\label{subs:excr}

We collect here some elementary facts about reversible Markov chains
on $\z_+$
which will later be used to give estimates for excursions of the
RWRE. Let $b \in \z_+$, $b>0$. Consider an excursion from $b$ to $b$.
Let $x \in \z_+$, $x>0$, $x \neq b$ and  denote by $Y_{b,x}$ the
number of visits to $x$ before returning to $b$.
The distribution of $Y_{b,x}$ is \lq\lq almost geometric\rq\rq : we have
$$
P_\omega(Y_{b,x} =m) =
\cases{ \alpha(1-\beta)^{m-1}\beta &$m=1,2,3, \ldots $, \cr\cr
        1-\alpha, &$m=0$, \cr}
$$
where $\alpha = \alpha_{b,x} = P_\omega(T(x) < T(b) \, | \, X_0 = b)$,
$\beta = \beta_{b,x} = P_\omega(T(b) < T(x) | X_0 = x)$. In
particular,
\begin{equation}
    E_\omega(Y_{b,x}) = {\alpha \over \beta} = {\mu(x) \over \mu(b)}
    = {\ee^{-V(x)} + \ee^{-V(x-1)} \over \ee^{-V(b)} + \ee^{-V(b-1)}}\,,
    \label{expu}
\end{equation}

\noindent where $\mu$ is the reversible measure for the Markov chain,
see (\ref{mu}). Further,
$$
\hbox{Var}_\omega(Y_{b,x}) = {\alpha (2- \beta - \alpha) \over
\beta^2}
\leq \frac{2}{\beta} {\mu(x) \over \mu(b)} .
$$ 
For $x >b
+1$, 
\begin{eqnarray*} 
\beta &=& (1-\omega_x)P_\omega(T(b)< T(x)| X_0 = x-1)\cr 
&=& (1-\omega_x)\Big(\sum\limits_{y=b}^{x-1} \ee^{V(y) - V(x-1)}\Big)^{-1} ,
\end{eqnarray*}
where the last formula follows from (\ref{hitbebeforei}),
and applies also for $x=b+1$. Hence,
for some $c_1=c_1(\delta)>0$, by (\ref{omega0}) and (\ref{mu}),
\begin{eqnarray}
    \hbox{Var}_\omega(Y_{b,x}) 
 &\le& c_1\, \ee^{-[\, V(x) - V(b)\,]}
\sum\limits_{y=b}^{x-1} \ee^{V(y) - V(x-1)}
    \nonumber
    \\
 &\le& c_1\, \ee^{-[\, V(x) - V(b)\,]} (x-b)\exp\Big( \max_{b\le
    y\le x-1} (V(y) - V(x-1))\Big) .   
    \label{varvis}
\end{eqnarray}
In the same way, one obtains, for $x < b$,
\begin{equation}
\label{varvis2}
\hbox{Var}_\omega(Y_{b,x}) \le c_1\, \ee^{-[\, V(x) - V(b)\,]} 
(b-x)\exp\Big( \max_{x \leq y\le b-1} (V(y) - V(x))\Big) . 
\end{equation}

\bigskip

\subsection{Number of valleys seen by the particle}
\label{subs:valley-seen}

Let $N_n$ denote 
the number of valleys ``seen" by the
particle in the first $n$ steps. More precisely, 
$$
N_n: = \sup\{k: \max_{0\le i\le n} X_i \geq m_k\}\,.
$$
Recall that as $k \to \infty$, 
$$
\frac{1}{2} \log m_k \sim \frac{1}{2} \log m_{k+1} \sim \log H_k \sim k
$$ 
(compare (\ref{mk1>Hk2}) with (\ref{bk<Hk2}) and use (\ref{logHk1})).
In combination with Fact \ref{f:revesz} this implies that
$\p$-almost surely,  
$$
\log \log T(m_k) \sim \log \log T(m_{k+1}) \sim \log H_k \sim k \,, 
\qquad k \to \infty\,.
$$
Since $T(m_{N_n}) \leq n < T(m_{N_n+1})$, it follows that 
\begin{equation}
\label{heightseen}
H_{N_n}  = (\log n)^{1+o(1)}\, , \qquad \hbox{\rm $\p$-a.s.} 
\end{equation}
and further 
\begin{equation}
    N_n \sim \log\log n, \qquad \hbox{\rm for } n\to \infty \qquad \hbox{\rm $\p$-a.s.} 
    \label{N}
\end{equation}

\subsection{The particle spends most of its time in the last two valleys}
\label{subs:spend}

Recall that $\xi(n,x)$ denotes the local time of the RWRE in $x$ at time $n$, and $m_k$ is the beginning of the $k$-th valley as in (\ref{mk}).
Let
\begin{equation}
    L(n, k) := \sum_{x\in [m_k,\; m_{k+1})} \xi(n,x) ,
    \label{L}
\end{equation}

\noindent which is the total time the particle spends
in the $k$-th valley during the first $n$ steps. 

The next theorem shows that the particle
spends most time in the two deepest valleys, which are the two right-most
valleys. 

\medskip

\begin{theorem}
 \label{t:spend} 
 We have, for any $\delta< 1$, 
 \begin{equation}
    \lim_{n\to \infty} \, {\exp\left( (\log n)^{\delta} \right)
    \over n} \, \sum_{1\le k< N_n-1} L(n,k) =0 \; , \qquad
    \hbox{\rm $\p$-a.s.}
    \label{negli}
 \end{equation}
 In particular,
 \begin{equation}
    \liminf_{n\to \infty} \, {1\over n} \sup_{k\ge 1} L(n,k)
    \ge {1\over 2} \; , \qquad \hbox{\rm $\p$-a.s.}
    \label{spend}
 \end{equation}
\end{theorem}

\medskip

\noindent {\it Proof.} It is clear that (\ref{spend}) follows from
(\ref{negli}) by taking $\delta=0$. 
Further, clearly (\ref{negli}) is a
consequence of (\ref{heightseen}) and 
\begin{equation}
    \lim_{N\to \infty} \, \ee^{(H_{N-2}^+)^\delta}\, \max_{n \in
    [T(m_N),\; T(m_{N+1}))} \, {1\over n} \,
    \sum_{1\le k< N-1} L(n,k) =0 \; , \qquad
    \hbox{\rm $\p$-a.s.}
    \label{tN}
\end{equation}
holding for any $\delta<1$.

In order to prove (\ref{tN}), we decompose the time interval
$[T(m_N),\; T(m_{N+1}))$ into excursions of the particle away from
$b_{N-1}$ and $m_{N-1}$.

Let $\varepsilon= \varepsilon_N >0$. Later, we will take
$\varepsilon_N= \exp(-(H_{N-2}^+)^\delta)$. Let
$$
n^* = n^*(N) := \inf \Big\{ n\ge T(m_N): \; \sum_{1\le k< N-1} L(n,k)
\ge \varepsilon \, n \Big\} ,
$$

\noindent with the notation $\inf\emptyset := \infty$. We are
interested in the case $n^* < T(m_{N+1})$; thus $n^*\in [T(i),
T(i+1))$ for some $i\in [m_N,\; m_{N+1})$. 

We define $T^1(b_{N-1}):= T(b_{N-1})$ and inductively,
\begin{eqnarray*}
    T^j(m_{N-1})
 &:=& \inf\left\{ n> T^j(b_{N-1}): \; X_n = m_{N-1}
    \right\},  
    \\
    T^{j+1}(b_{N-1})
 &:=& \inf\left\{ n> T^j(m_{N-1}): \; X_n = b_{N-1}
    \right\}, \qquad j\ge 1.
\end{eqnarray*}

\noindent For any $i\in [m_N,\; m_{N+1})$, let $M_i := \sup\{ j: \;
T^j(m_{N-1}) < T(i+1)\}$ (notation: $\sup \emptyset := 0$), be the total
number of excursions from $b_{N-1}$ to $m_{N-1}$, before reaching
$i+1$. 

If $n^*\in [T(i), T(i+1))$ and $M_i=0$, we 
have
$\sum_{1\le k< N-1}
L(n^*,k) \le T(b_{N-1})$ and $n^* \ge T(i) - T(b_{N-1})$ so that
\begin{equation}
    T(b_{N-1}) \ge \varepsilon (T(i) - T(b_{N-1})) ;
    \label{m=0}
\end{equation}

\noindent whereas if $n^*\in [T(i), T(i+1))$ and $M_i\ge 1$, then
$\sum_{1\le k< N-1} L(n^*,k) \le T^1(b_{N-1}) + \sum_{j=1}^{M_i}
[ \, T^{j+1}(b_{N-1}) - T^j(m_{N-1}) \, ]$ and $n^* \ge 
\sum_{j=1}^{M_i} [ \, T^j(m_{N-1}) - T^j(b_{N-1})\, ]$ so that
\begin{equation}
    T^1(b_{N-1}) + \sum_{j=1}^{M_i} [ \, T^{j+1}(b_{N-1}) - T^j(m_{N-1}) \,
    ] \ge \varepsilon \sum_{j=1}^{M_i} \left[ T^j(m_{N-1}) -
    T^j(b_{N-1}) \right].
    \label{m>=1}
\end{equation}

We first treat the case $M_i = 0$, i.e., there is no excursion (before
time $T(i+1)$) back to $m_{N-1}$ after reaching $b_{N-1}$. In this
case, (\ref{m=0}) holds. Let
\begin{eqnarray}
    p_{i,N}
 &:=& P_\omega\left(T(b_{N-1}) \ge \varepsilon (T(i) -
    T(b_{N-1}))\right)
    \nonumber \\
 &=& P_\omega \Big( T(b_{N-1}) \ge {\varepsilon \over 1+ \varepsilon} 
    T(i) \Big)
    \nonumber
    \\
 &\le& P_\omega \left( T(b_{N-1}) \ge \lambda\,\right)
    +P_\omega \Big( T (i) < {(1+\varepsilon)\lambda \over
    \varepsilon}\Big) ,
    \label{pnest}
\end{eqnarray}

\noindent for any $\lambda\ge 1$. Considering the first term in
(\ref{pnest}), we have, by (\ref{cheby}),
$$
P_\omega \left(T(b_{N-1}) \ge \lambda\,\right) \le {b_{N-1}^2 \over
\lambda}\exp(H_{N-2}^+) .
$$
Turning to 
the second term in (\ref{pnest}), we have 
\begin{eqnarray*}
P_\omega \Big( T (i) < {(1+\varepsilon)\lambda\over 
\varepsilon}\Big)
&\le&
P_\omega \Big( T (m_{N})< {(1+\varepsilon)\lambda\over 
\varepsilon}\; \Big| \; X_0=b_{N-1} \Big)\\
&\le&
P_\omega \Big( \tau_{N-1} < {(1+\varepsilon)\lambda\over 
\varepsilon}\, \Big)\\
&\le& {
c_0 
(1+\varepsilon)\lambda\over 
\varepsilon} \ee^{- H_{N-1}},
\end{eqnarray*}
where we used (\ref{exit}) for the last inequality. 
Hence, plugging in the value of $\varepsilon =
\ee^{-(H_{N-2}^+)^\delta}$,
\begin{eqnarray*}
p_{i,N} 
&\le&
{b_{N-1}^2 \over \lambda}\exp(H_{N-2}^+)
+ {c_0 (1+\varepsilon)\lambda\over 
\varepsilon}\exp(- H_{N-1})\\
&\le&
{b_{N-1}^2 \over \lambda}\exp(H_{N-2}^+)
+2 c_0\lambda \exp((H_{N-2}^+)^\delta)\exp(- H_{N-1}) .
\end{eqnarray*}
We choose $\lambda =\lambda_N:= \exp(\hetzi H_{N-1}+\hetzi
H_{N-2}^+)$. Then, 
$$
p_{i,N} \leq (b_{N-1}^2 +2 c_0) \exp\left(-\hetzi H_{N-1}+\hetzi
H_{N-2}^+ + (H_{N-2}^+)^\delta\right)\, .
$$

\noindent Due to (\ref{bk<Hk2}) and Lemma \ref{l:Hk}, $b_{N-1} \le
(H_{N-1}^+)^{3}$ and $m_{N+1} \le (H_{N-1})^3$ for $N\to \infty$, so
that by Lemmas \ref{l:Hk-Hk-1} and \ref{l:Hk}, 
\begin{equation}
\label{zeroexc}
\sum_N\sum_{m_N \le i < m_{N+1}} P_\omega\left( n^*\in [T(i),   
T(i+1)), \; M_i=0\right) <\infty, \quad
\hbox{$P$-a.s.} 
\end{equation}

Turning to consider $n^*\in [T(i), T(i+1))$ and $M_i \ge 1$, for 
$\lambda=\lambda_N>0$ to be chosen later, and each $m \geq 1$ let
\begin{eqnarray*} 
A(m) &:=&  \Big\{ T^1(b_{N-1})+ \sum_{j=1}^{m} \left[ T^{j+1}(b_{N-1}) -
    T^j(m_{N-1}) \right] \ge m\lambda \Big\} \\
B(m) &:=&   \Big\{ 
    \sum_{j=1}^m \left[ T^j(m_{N-1}) - T^j(b_{N-1}) \right] <
    {m\lambda\over \varepsilon} \Big\} \,. 
\end{eqnarray*}    
Note that if $n^*\in [T(i), T(i+1))$ for some 
$i \geq m_N$ with $M_i \ge 1$, then (\ref{m>=1}) holds,
and hence either $A(M_i)$ or $B(M_i)$ holds as well. 
Consequently, decomposing the event $A(M_i)$ according to $i$
and the event $B(M_i)$ according to the value $m$ of $M_i$, we get that
\begin{eqnarray}
 &&P_\omega\left( n^* \in [T(i), \; T(i+1)), \hbox{ \rm for some }i
    \in [m_N,\; m_{N+1})\hbox{ \rm and  } M_i \geq 1\right)
    \nonumber \\
 &\le& \sum_{i=m_N}^{m_{N+1}} P_\omega \left(A(M_i), M_i \geq 1\right) +
  P_\omega \left(B(M_j) \hbox{ \rm for some } j \geq  m_N
\hbox{ \rm and  } M_j \geq 1\right)
\nonumber \\
&\le& m_{N+1} \sup_{i \geq b_{N-1}}
P_\omega\left( A(M_i), M_i \geq 1 \right) 
+ \sum_{m=1}^{\infty} P_\omega\left(B(m)\right) \nonumber \\ 
 &=:& 
m_{N+1} I^{(1)} + \sum_{m=1}^\infty I_m^{(2)}.
    \label{P(L>L)}
\end{eqnarray}
By the strong Markov property, conditionally on $\omega$ 
both $T^1(b_{N-1})$ and the identically distributed 
random variables 
$T^{j+1}(b_{N-1}) - T^j(m_{N-1})$, $j\ge 1$,
are independent of the value of $M_i$ for $i \geq b_{N-1}$. Hence,
by Markov's inequality
\begin{eqnarray}
I^{(1)} 
&\leq& \sup_{m \geq 1} P_\omega \left(A(m)\right) 
\le \sup_{m \geq 1} 
{1 \over m\lambda}E_\omega\left(T(b_{N-1}) \right)+
 {1 \over \lambda} E_\omega\left(T(b_{N-1}) \, | \, X_0=
    m_{N-1} \right) 
\nonumber \\
 &\le& {2\over \lambda} E_\omega\left(T(b_{N-1}) \right) 
 \le {2b_{N-1}^2\over \lambda} \exp\left( H_{N-2}^+ \right) ,
\label{I1bd}
\end{eqnarray}
where the last inequality is due to 
(\ref{l:hitting}). 

Further, since
$T^j(m_{N-1}) - T^j(b_{N-1})$, $j\ge 1$,
are i.i.d.\ random
variables, each having the law of $T(m_{N-1})$ when starting 
at $b_{N-1}$, by Corollary \ref{c:exit}, for any $a>0$, 
$$
I_m^{(2)} 
 \le \ee^{am\lambda/\varepsilon} \left( E_\omega \Big(
\ee^{- a T (m_{N-1})} \; \big| \; X_0 = b_{N-1} \Big) \right) ^{\! m}
\le \Big( \ee^{a\lambda/\varepsilon} {
2 c_0 \ee^{- H_{N-1}}\over a}\Big)^{\!\! m}.
$$
We choose $\lambda=\lambda_N:=  \exp\left(\, H_{N-1} -
2(H_{N-2}^+)^\delta\, \right)$, $a=a_N:= \exp \left(\, -H_{N-1} +
(H_{N-2}^+)^\delta\, \right)$ and as stated 
before $\varepsilon = \varepsilon_N
:= \exp \left(\, -(H_{N-2}^+)^\delta\, \right)$. Since
$a \lambda \varepsilon^{-1} =1$ these choices result with 
\begin{equation}\label{eq:amir-I1}
\sum_{m=1}^\infty I_m^{(2)} \le 
\sum_{m=1}^\infty
 \Big( 2 \ee c_0 \exp(- (H_{N-2}^+)^\delta) \Big)^{\!\! m} \leq c_2 
 \exp( - (H_{N-2}^+)^\delta) \;. 
\end{equation}
In view of (\ref{I1bd}), these choices also lead to 
$$
m_{N+1} I^{(1)} \leq  
2 m_{N+1} b_{N-1}^2 \exp( 
-H_{N-1} +2(H_{N-2}^+)^\delta + H_{N-2}^+) \,,
$$
where $P$-almost surely, for all large $N$
$$
-H_{N-1} + 2 (H_{N-2}^+)^\delta +
H_{N-2}^+ \le - (H_{N-2}^+)^{2\delta} 
$$  
(see Lemma \ref{l:Hk-Hk-1}). Further, 
due to (\ref{bk<Hk2}) and Lemma \ref{l:Hk}, 
$P$-almost surely, for all large $N$,
$$
2 m_{N+1} b_{N-1}^2 \le (H_{N-2}^+)^7 \,,
$$
yielding that  
$$
m_{N+1} I^{(1)} \leq \exp(-(H_{N-2}^+)^\delta) \,. 
$$
Combining this with (\ref{eq:amir-I1})
and (\ref{P(L>L)}) yields, together with
(\ref{zeroexc}), that
$$
\sum_N P_\omega \left( n^* \in [T(m_N), \; T(m_{N+1})) \right)  
<\infty .
$$
By the Borel--Cantelli lemma, we obtain that for any $\delta<1$,
$$
\limsup_{N\to \infty} \, \ee^{(H_{N-2}^+)^\delta}\, \max_{n \in
[T(m_N),\; T(m_{N+1}))} \, {1\over n} \, \sum_{1\le k< N-1} L(n,k) \le
1 \; ,
\qquad \hbox{\rm $\p$-a.s.}
$$
Since $\delta<1$ is arbitrary, this implies (\ref{tN}),
and completes the proof of Theorem \ref{t:spend}. \hfill$\Box$

\section{Occupation time and local time}
\label{s:lt}

We have so far proved in Theorem \ref{t:spend} that ($\p$-almost
surely for $n$ large enough) the particle spends at least $({1\over
2}+o(1))n$ time in a certain valley. The goal of this section is
to prove that the time spent by the particle at the bottom of this or 
a neighbor valley is at least a constant multiple of $n/\log \log\log
n$.

There are two main points in the proof: (a) We need to
investigate the ratio between the time spent in a valley
(occupation time) and the time spent in the bottom of the same (or a 
neighbor) valley (local time). This is the main part of this section;
(b) Since the valley where the particle spends at least $({1\over
2}+o(1))n$ time has a random number (namely, $N_n$ or $N_n-1$, see
Section \ref{s:particle}) and this random number depends on the
environment as well as on the movement of the particle, we need a
result which holds uniformly for a whole collection of valleys.

\subsection{Comparison between occupation time and local time}
\label{subs:ratio}

Recall that $N_n$ is the number of valleys seen by the particle in
the first $n$ steps. Define, for any $k\ge 1$,
\begin{equation}
    \Lambda_k:= \sum_{i=m_k}^{m_{k+1}-1} \ee^{-[\, V(i)-V(b_k)\,]} \; .
    \label{Lambda}
\end{equation}

\noindent Note that $(\Lambda_k, \; k\ge 1)$ depends only on the
environment, and that
\begin{equation}
    \inf_{k\ge 1} \Lambda_k \ge 1.
    \label{Lambda>c}
\end{equation}

\noindent Here are the main estimates of this subsection, which relate
occupation time with local time.
In particular note that $\Lambda_k$ measures the effective width of 
the $k$-th valley as reflected by the ratio between the 
expected
occupation time and the maximal 
expected
local time among its sites
(at the appropriate time $n=T(m_{k+1})$ of the particle just 
reaching the beginning of the next valley).

\medskip

\begin{proposition}
 \label{l:ratio} 
 There exist $c_3$ and $c_4$ such that $\p$-almost
 surely for $n$ large enough,
 \begin{eqnarray}
     L(n, N_n-1) 
  &\le& c_3 \, \Lambda_{N_n-1} \, \xi(n, b_{N_n-1}),
     \label{ratio}
     \\
     L(n, N_n) 
  &\le& c_4 \, \Lambda_{N_n} \, \left[\, \xi(n, b_{N_n-1}) 
     + \xi(n, b_{N_n}) \, \right],    
     \label{ratioN}
 \end{eqnarray}
 where $L(n,k)$ is the time spent in the $k$-th valley as in
 $(\ref{L})$.
\end{proposition}

\begin{proposition}
 \label{l:ratio2} 
  There exists $c_5$ such that $\p$-almost surely for all
  large $N$,
 \begin{equation}
    L(T(m_N), N-1) \ge c_5 \, \Lambda_{N-1} \, \max_{x\in
    [m_{N-1}, m_N)} \xi(T(m_N), x).
    \label{ratio-lb}
 \end{equation}
\end{proposition}

\medskip

\noindent {\it Remark on the proof.} The basic idea of the proof of
the propositions can be described as follows. For (\ref{ratio-lb}), we
consider excursions of the walk away from $b_{N-1}$ during the time
interval $[T(b_{N-1}), T(m_N)]$, and let $M=M(N)$ denote the number
of completed excursions ($M$ can be 0). The random variable
$M$, which is $\xi(T(m_N), b_{N-1}) -1$, has a geometric distribution
(under $P_\omega$)
and $E_\omega(M)$ is approximatively
$\ee^{H_{N-1}^+}$.
 By the strong Markov property, all completed
excursions make i.i.d. contributions to $\xi(T(m_N), x)$, for any $x$,
hence also to $L(T(m_N), N-1)$.  The law of large numbers says that,
with $\rho$ denoting the lifetime of an excursion,
$$
    \xi(T(m_N), x) \asymp M\, E_\omega(\xi(\rho,x)) \asymp
    \xi(T(m_N), b_{N-1})\ee^{-[V(x) - V(b_{N-1})]},
$$

\noindent (it was proved in Subsection \ref{subs:excr} that
$E_\omega(\xi(\rho,x)) \asymp \ee^{-[V(x) - V(b_{N-1})]}$), and
similarly, 
\begin{equation}
    L(T(m_N), N-1) \asymp M \sum_{x\in [m_{N-1}, m_N)} E_\omega(\,
    \xi(\rho,x)\,)  \asymp \xi(T(m_N), b_{N-1}) \, \Lambda_{N-1}.
    \label{L=xi(b)}
\end{equation}

\noindent This would yield (\ref{ratio-lb}) if we take $c_5$ to be
sufficiently small. In order to give a rigorous proof of
(\ref{ratio-lb}), we need to estimate deviation probabilities for $M$
(which is easy), and for the number of visits during a single
excursion (which is done via a second moment argument).

\medskip
The proof of Proposition \ref{l:ratio} needs slightly more care since it
involves an arbitrary time $n$, instead of the first hitting times
$T(m_N)$ in Proposition \ref{l:ratio2}. Both proofs go
along the lines described in the preceding remark, but require 
certain technical adjustments. We start with a few preliminary estimates. 
The first is a rigorous statement of (\ref{L=xi(b)}). For further needs
we now provide such a statement uniformly over all $n \geq T(m_N)$, 
instead of just for $T(m_N)$.

\begin{lemma}
\label{l:L=xi(b)}
There exist $0<c_6<c_3<\infty$ such that, for any $\varepsilon>0$, 
$P$-almost surely for all $N$ large enough, 
\begin{eqnarray}
P_\omega\left(\exists n \ge T(\eta_N), \;
L(n, N) \ge c_3 \, \Lambda_N \xi(n, b_N) \right)
&\le& \ee^{- (H_N)^{1-\varepsilon}},
\label{L>xi(b)}
\\
P_\omega\left(\exists n \geq T(m_{N+1}), \;
L(n,N) \le c_6 \, \Lambda_N\xi(n,b_N) \right)
&\le& \ee^{- (H_N)^{1-\varepsilon}}.
\label{L<xi(b)}
\end{eqnarray}
\end{lemma}

\medskip

\noindent {\it Proof of Lemma \ref{l:L=xi(b)}.}
We decompose the random walk into excursions
away from $b=b_N$. That is, $T^{-1}=0$, $T^0:= T(b)$ and
\begin{equation}
    T^j 
 :=\inf \left\{ k> T^{j-1}: \; X_k = b \right\},
    \label{T}
\end{equation}
are the times of consecutive visits to $b$,
which are 
$\p$-almost
surely finite on 
account of (\ref{recurrent}). Fixing $i \in [\eta_N,m_{N+1})$, 
consider the corresponding occupation times of the interval 
$[m_N,i]$, that is, 
$$
Z_j = Z_j (i) := \sum_{x=m_N}^i \xi(T^j,x)-\xi(T^{j-1},x) \,.
$$
Note that, by the strong Markov property of the walk,
$Z_j$, $j\ge 1$,
 are 
independent non-negative random variables
(under $P_\omega$), 
and are
 also identically 
distributed and of finite second moment (c.f.
(\ref{eq:bd-var}) in the sequel). Observe that
$$
\oo{M}:= \xi(n,b) = \inf \left\{ j: \; T^j > n \right\},
$$ 
and $\oo{M} \geq 1$ whenever $n \geq T(b)$ (which is always the 
case here). Further, for $i=m_{N+1}-1$,
\begin{equation}\label{eq:yuval-bd}
L(n,N) \leq Z_0 + \sum_{j=1}^{\oo{M}} Z_j \,,
\end{equation}
and (\ref{eq:yuval-bd}) applies also for $i<m_{N+1}-1$ provided $n<T(i+1)$.

Since $Z_j \geq 0$, it follows that for any $i \in [\eta_N,m_{N+1})$,
$c_7>0$, $\ell \geq 1$ and $k_r = \ell 2^r$,
\begin{eqnarray}
\nonumber 
&& P_\omega\big(\exists n \in [T(i),T(i+1)), \;
L(n, N) \ge (2c_7+1) \oo{M} \Lambda_{N} \big)
\\
&\le& P_\omega (\exists n \geq T(i), \;
\oo{M}  \leq \ell) + P_\omega \big(Z_0 \ge \ell \Lambda_N \big)
+ \sum_{r=0}^\infty
P_\omega \big(\sum_{j=1}^{k_r} Z_j \ge c_7 k_r \Lambda_N \big) 
\label{M>ell}
\\
&=:&
 I_1(i) + I_2 + I_3(i) \,.
\nonumber 
\end{eqnarray}
Further, as the inequality (\ref{M>ell})
holds for $i=m_{N+1}-1$ even without the condition $n<T(i+1)$, we have 
for $c_3=2c_7+1$ that 
\begin{equation}\label{M>ell-rev}
P_\omega\Big(\exists n \geq T(\eta_N), \;
L(n, N) \ge c_3 \oo{M} \Lambda_{N} \Big) \le
\sum_{i=\eta_N}^{m_{N+1}-1} (I_1(i) + I_2 + I_3(i)) \,.
\end{equation}

To estimate the term $I_1(i)$ in (\ref{M>ell}), let $K(b,i)$ denote the number 
of excursions from $b$ to $b$
made by the walk during the time interval $[T(b),T(i)]$, 
which has a
geometric distribution of parameter $p=p(b,i)$, that is, $P_\omega
(K=k) = (1-p)^k p$, $k=0,1,2, \ldots$, where due to
(\ref{hitbebeforei}), for any $i>b$, 
\begin{equation} p(b,i) :=
\omega_{b} P_\omega ( T(i)<T(b) \, | \, X_0 = b+1) = \omega_{b} \,
{\ee^{V(b)} \over \sum_{y=b}^{i-1} \ee^{V(y)}} \leq \ee^{-W(b,i)} \,,
\label{geom1} 
\end{equation} 
for $W(b,i):=\max_{b \leq y <i}
V(y)-V(b)$. In particular, for $i \geq \eta_N$ we have that
\begin{equation} 
p(b_N,i) \leq \ee^{-W(b_N,\eta_N)} \leq
\ee^{V(b_N)-V(\eta_N) + C} \leq c_8 \ee^{-H^+_{N-1}} \,.  
\label{geom}
\end{equation} 
$~$ For any $i>b$, the event $\{ n \geq T(i) \}$ implies that
$\oo{M} > K(b,i)$. Hence, fixing $\varepsilon>0$ and
$\ell:= \lceil p(b,i)^{-1} \exp(-\frac{1}{3}(H_{N-1}^+)^{1-\varepsilon})
\rceil$, we have that 
\begin{equation} \label{4-1} 
I_1(i) \leq P_\omega
(K(b,i) < \ell) = 1-(1-p)^{\ell} \le p\, \ell \le 
c_8 \ee^{-\frac{1}{3} (H_{N-1}^+)^{1-\varepsilon}} \, .  
\end{equation}

Proceeding to deal with $I_2$, 
since the steps of the random walk within $[0,m_N-1]$ do not matter
to $Z_0(i)=Z_0(b_N)$, the latter has under $P_\omega$ the same 
law as that of the occupation time of $[1,b_N-m_N+1]$ 
till $T(b_N-m_N+1)$ under $P_{\widetilde{\omega}}$, where 
$\widetilde{\omega}_x = \omega_{x+m_N-1}$. Consequently,
by (\ref{l:hitting}) we have that $P$-almost surely, for all $N$ large
enough,
$$
E_\omega (Z_0) \le E_{\widetilde{\omega}} (T(b_N-m_N+1)) \leq
b_N^2 \exp\big( \, \max_{m_N-1 \le
y\le z< b_N} (V(z)-V(y)) \big) \le 
b_N^2\, \ee^{H_{N-1}^+ - (H_{N-1}^+)^{1-\varepsilon}},
$$
with the last inequality due to (\ref{depth1}).  
It follows that for our choice of $\ell=\ell(i,N,\varepsilon)$,
\begin{equation}
I_2 \leq P_\omega( Z_0 \geq \ell) \le  \ell^{-1} E_\omega (Z_0) \le
c_8 b_N^2 \ee^{- \frac{2}{3} (H_{N-1}^+)^{1-\varepsilon}} 
\label{4-2}
\end{equation}
(where the first inequality is due to (\ref{Lambda>c}) and the last
one due to (\ref{geom})).

As for the term $I_3(i)$ of (\ref{M>ell}),
observe that in the notations of Subsection \ref{subs:excr},
$$
Z_1 = \sum_{x=m_N}^i [\, \xi(T^1,x) - \xi(T^0, x)\, ] 
= \sum_{x=m_N}^i Y_{b,x} \,,
$$
where, by (\ref{expu}), 
\begin{equation}
\label{constbound}
E_\omega\left(\, \xi(T^1, x) - \xi(T^0,x)\, \right) = E_\omega\left(Y_{b,x}\right) 
= {\omega_b\over \omega_x} \ee^{-[V(x)-V(b)]} \, .
\end{equation}
It follows, in view of assumption (\ref{omega0}), that
\begin{equation}
E_\omega (Z_1) \le \delta^{-1} \, \Lambda_N.
    \label{cst}
\end{equation}
Consequently, by the independence of $Z_j$ we get
for $c_7 \geq \delta^{-1}+1$ and $k_r=\ell 2^r$, the bound
\begin{equation}\label{eq:I3}
I_3(i) \leq \sum_{r=0}^\infty
P_\omega \big(\sum_{j=1}^{k_r} (Z_j-E_\omega (Z_j)) \ge
 k_r \Lambda_N \big)
\leq \frac{\hbox{Var}_\omega (Z_1)}{\Lambda_N^2} \sum_{r=0}^\infty 
\frac{1}{k_r}
\leq \frac{2 \hbox{Var}_\omega (Z_1)}{\ell} \,,
\end{equation}
using (\ref{Lambda>c}) in the last inequality.
Observe that 
\begin{eqnarray*}
 \hbox{Var}_\omega (Z_1)
 = \hbox{Var}_\omega \big( \sum_{x=m_N}^i Y_{b,x} \big)
 \le m_{N+1} \sum_{x=m_N}^i \hbox{Var}_\omega ( Y_{b,x} ) \,.
\end{eqnarray*}
Since $V(x) \geq V(b)$ for all $x\in [m_N,m_{N+1})$,
we have from (\ref{varvis}) that for $b=b_N$ and any $x\in (b, \, m_{N+1})$,
\begin{equation}
 \hbox{Var}_\omega (Y_{b,x})
 \le c_1 m_{N+1} \exp\Big(
    \max_{b\le y<x} (V(y)-V(x-1))\Big) \,.
    \label{var-excr}
\end{equation}
Similarly, applying (\ref{varvis2}) instead of (\ref{varvis}), 
we obtain that for all $x\in [m_N, \, b)$,
\begin{equation}
    \hbox{Var}_\omega (Y_{b,x}) \le
     c_1 \, m_{N+1} \exp\Big(
    \max_{x \le z<b} (V(z)-V(x))\Big) \,,
    \label{var-excr3}
\end{equation}
and of course $\hbox{Var}_\omega (Y_{b,b}) = 0$. 
Summing over $x \in [m_N,i]$ we find
by means of (\ref{var-excr}) and (\ref{var-excr3}) that 
\begin{equation}\label{eq:bd-var}
\hbox{Var}_\omega (Z_1(i)) \le c_{9}\, m_{N+1}^3 \ee^{U(b,i)} \,,
\end{equation}
where 
$$
U(b,i):=\max\{ \max_{m_N \le y \le z<b} (V(z)-V(y)), \; 
\max_{b \le y\le z < i} (V(y)- V(z)) \}\,.
$$
Let
$$
\Delta_N := \max_{i \in [\eta_N,m_{N+1})} \{ U(b,i)-W(b,i) \} \,.
$$
Combining (\ref{eq:I3}) and (\ref{eq:bd-var}) we see that
by (\ref{geom1}), for our choice of $\ell=\ell(i,N,\varepsilon)$,
\begin{equation}\label{4-3}
I_3(i) \leq 2 p(b_N,i) \hbox{Var}_\omega (Z_1(i))
\ee^{\frac{1}{3}(H_{N-1}^+)^{1-\varepsilon}}
\leq c_{10} m_{N+1}^3 \ee^{\Delta_N +
\frac{1}{3}(H_{N-1}^+)^{1-\varepsilon}} \,.
\end{equation}
Combining (\ref{depth1}) and (\ref{depth4}), we deduce that
$P$-almost surely, for all $N$ large enough,
$$
U(b,\eta_N) \le H_{N-1}^+ - (H_{N-1}^+)^{1-\varepsilon} \,.
$$
Likewise, note that if $i\in (\eta_N, m_{N+1})$ then combining
the preceding with (\ref{depth3}) we have that
$$
U(b,i) \le \max \{ U(b,\eta_N),
W(b,i) - (H_{N-1}^+)^{1-\varepsilon} \} \le
W(b,i) - (H_{N-1}^+)^{1-\varepsilon} \,,
$$
using in the last inequality the fact that if $i > \eta_N$ then
$$
W(b,i) = \max_{b \leq y < i} V(y) - V(b) \geq V(\eta_N)-V(b) \geq H_{N-1}^+ \,.
$$
Consequently, $P$-almost surely, for all $N$ large enough
$$
\Delta_N \leq C - (H_{N-1}^+)^{1-\varepsilon} \,,
$$
and thus,
plugging (\ref{4-1}), (\ref{4-2}) and (\ref{4-3}) into (\ref{M>ell-rev}),
we obtain that $P$-almost surely, for all $N$ large 
$$
P_\omega\Big(\exists n \geq T(\eta_N), \;
L(n, N) \ge c_3 \oo{M} \Lambda_{N} \Big) \le
c_{11} m_{N+1}^4 \ee^{-\frac{1}{3} (H_{N-1}^+)^{1-\varepsilon}} \, .  
$$
Noting that $\varepsilon>0$ is arbitrary, in view of (\ref{bk<Hk2}) and
Lemma \ref{l:Hk}, this implies (\ref{L>xi(b)}).

Moving next to the proof of (\ref{L<xi(b)}), since we are not 
considering $n < T(m_{N+1})$ in this inequality, we set 
$i=m_{N+1}-1$ for the remainder of the proof, in which case 
we have from (\ref{constbound}) that 
$$
E_\omega(Z_1) =\sum_{x=m_N}^{m_{N+1}-1} {\omega_{b_N}\over \omega_x} 
\ee^{-[V(x)-V(b_N)]} \ge \delta \sum_{x=m_N}^{m_{N+1}-1} 
\ee^{-[V(x)-V(b_N)]}  \geq \delta \Lambda_N \,.
$$ 
Since further $L(n,N) \geq \sum_{j=1}^{\oo{M}-1} Z_j$ in this case
(regardless of $n$), we have similarly to (\ref{M>ell}),
that for any $c_6>0$ and $k_r=\ell 2^r$, $\ell\ge 1$,
\begin{eqnarray}
&& P_\omega \Big(\exists n \geq T(m_{N+1}), \; L(n, N) 
\le c_6 \oo{M} \Lambda_{N} \Big) \nonumber \\
 &\le& P_\omega (\exists n \geq T(m_{N+1}), \oo{M} \leq \ell)
+ \sum_{r=0}^\infty P_\omega
    \label{M<ell}
\big(\sum_{j=1}^{k_r} Z_j \le 4 c_6 k_r \Lambda_N \big) \\
&=:&
  I_1 + I_4 \,.
\nonumber 
\end{eqnarray}
With $E_\omega(Z_1) \geq \delta \Lambda_N$, note that
if $c_6 < \delta/5$, then
\begin{eqnarray}
I_4 &\leq&
 \sum_{r=0}^\infty
P_\omega
\Big(\sum_{j=1}^{k_r}
(Z_j-E_\omega (Z_j)) \le - \frac{\delta}{5} k_r \Lambda_N \Big)
\le {c_{12} \hbox{Var}_\omega (Z_1) \over \ell} \,,
\label{Var(eta)}
\end{eqnarray}
(using in the last inequality both (\ref{Lambda>c}) and
the fact that $\sum_r k_r^{-1} = 2 \ell^{-1}$). Thus,
taking $\ell=\ell(i,N,\varepsilon)$ as before, in view of
(\ref{M<ell}) and (\ref{Var(eta)}) we get (\ref{L<xi(b)})
by the same argument used to complete the derivation of (\ref{L>xi(b)})
(even simpler, as we neither sum over $i$ nor consider $I_2$ here).
\hfill$\Box$

\medskip
We next show that upon the walk reaching 
the right end of a given valley, with high probability 
no point of this valley has a local time much 
larger than its bottom. This estimate 
complements (\ref{L<xi(b)}) en-route to
proving Proposition \ref{l:ratio2}.

\begin{lemma}
\label{l:xi>xi*} 
There exists $\gamma$ finite 
such that for any $\varepsilon>0$,  
$P$-almost surely for $N$ large enough,
$$ 
\max_{x\in [m_N, m_{N+1})}
P_\omega \Big( \xi(T(m_{N+1}), x) \ge \gamma \xi(T(m_{N+1}), b_N) \Big) \le
    \ee^{-(H_N)^{1-\varepsilon}}. 
$$ 
\end{lemma}

\medskip
\noindent {\it Proof of Proposition \ref{l:ratio2}.} Taking 
$c_5=c_6/\gamma>0$, the proposition follows from 
(\ref{L<xi(b)}) and Lemma \ref{l:xi>xi*} by means of the
Borel--Cantelli lemma (as 
$m_N \ee^{-(H_N)^{1-\varepsilon}}$ is summable). 
\hfill$\Box$ 

\medskip
\noindent {\it Proof of Lemma \ref{l:xi>xi*}.} We use the
path decomposition of the walk as in Lemma \ref{l:L=xi(b)},
with $b=b_N$, $T^{-1}=0$, $T^0=T(b)$ and $T^j$, $j \geq 1$
the times of returns of the walk to $b$ (c.f. (\ref{T})). We
further set $\oo{M}=\xi(n,b) \geq 1$, hereafter taking $n=T(i+1)$ for
$i=m_{N+1}-1$. Fixing $x\in [m_N,i]$, $x \neq b$,
let $Y_j = \xi(T^j,x) - \xi(T^{j-1}, x)$ for $j=0,1,\ldots$,
denote the accumulated local time at $x$ during 
the $j$-th segment of the walk. Note 
that the non-negative random variables $Y_j$, $j \geq 1$, are i.i.d. 
of finite second moment, and with
$Y_1$ having the law of $Y_{b,x}$ of Subsection \ref{subs:excr}, 
also $E_\omega (Y_1) \leq \delta^{-1}$ (c.f. (\ref{constbound})). Further,
similarly to (\ref{eq:yuval-bd}) we have that 
$$
\xi(n,x) \leq Y_0 + \sum_{j=1}^{\oo{M}} Y_j \,.
$$
Hence, as in (\ref{M>ell}), for $n=T(i+1)$,
$\gamma \geq 2(1+\delta^{-1})+1$ and $k_r = \ell 2^r$, $\ell \geq 1$,
\begin{eqnarray*}
P_\omega (\xi(n,x) \ge \gamma \oo{M} )
 &\le& P_\omega (\exists 
n \geq T(i), \oo{M} \leq \ell) + P_\omega (Y_0 \geq \ell)
+ \sum_{r=0}^\infty
P_\omega \Big(\sum_{j=1}^{k_r} (Y_j - E_\omega (Y_j)) \geq k_r\Big)
\\
&=:&
 I_1 + I_2 + I_5(x) \,.
\end{eqnarray*}
We fix $\varepsilon >0$ and
$\ell=\ell(i,N,\varepsilon)$ as in Lemma \ref{l:L=xi(b)},
thus taking care of the term $I_1$ (c.f. (\ref{4-1})). Further,
with $Y_0 \leq Z_0$ this choice also takes care
of $I_2$ (c.f. (\ref{4-2})) and just as in
(\ref{eq:I3}) we have that
$$
I_5(x) \leq \frac{2 \hbox{Var}_\omega (Y_1)}{\ell} =
 \frac{2 \hbox{Var}_\omega (Y_{b,x})}{\ell} \,.
$$
It follows from (\ref{var-excr}) and (\ref{var-excr3}) that
$$
\max_{x \in [m_N,i]} \, \hbox{Var}_\omega (Y_{b,x}) \leq
c_9 m_{N+1} \ee^{U(b,i)}
$$
(compare with the derivation of (\ref{eq:bd-var})). For 
our choice of $\ell$ and the bound (\ref{geom1})
on $p(b,i)$ it follows that
$P$-almost surely, for any $N$ large enough and 
all $x \in [m_N,i]=[m_N,m_{N+1})$,
$$
I_5(x) \leq 
c_{13} m_{N+1} \ee^{-\frac{2}{3} (H_{N-1}^+)^{1-\varepsilon}}
$$
(see (\ref{4-3})). As observed before, such estimates
are all we need for the lemma (in view of (\ref{bk<Hk2}) 
and Lemma \ref{l:Hk}).
\hfill$\Box$

Our next lemma is similar in spirit to Lemma \ref{l:L=xi(b)}. Its proof 
is slightly more involved since two different (consecutive) valley 
bottoms are relevant here. This happens for example when the occupation 
time of the last seen valley is to be considered, as in (\ref{ratioN}).

\medskip

\begin{lemma}
 \label{l:xi+xi} 
 There exists $\kappa$ finite such that for any $\varepsilon>0$, 
 $P$-almost surely for $N$ large enough,
 $$
    P_\omega \Big(\exists n < T(\eta_N), \;
    L(n, N) > \kappa \, [\, \xi(T(m_N),
    b_{N-1}) + \xi(n, b_N)\, ] \Lambda_N \Big)
    \le \ee^{-(H_N)^{1-\varepsilon}}.  
 $$ 
\end{lemma}

\medskip

\noindent {\it Proof of Lemma \ref{l:xi+xi}.}
Clearly, it suffices to consider $n \in [T(i),T(i+1))$ for
$i \in [m_N,\eta_N)$. To this end, we 
adopt the path decomposition and notations of Lemma \ref{l:L=xi(b)}
(for $b=b_N$). The random variables
$Z_j$, $j \geq 1$ are i.i.d. and the basic
inequality (\ref{eq:yuval-bd}) applies, just taking
$\oo{M}=0$ whenever $n \in [T(i),T(i+1))$ for $i <b$.
With $\oo{K} :=\xi(T(m_N), b_{N-1})$, recall that $K=\oo{K}-1$ is a
geometric random variable of parameter 
\begin{equation}\label{eq:p-bd}
p(b_{N-1},m_N) \leq \ee^{-W(b_{N-1},m_N)} \leq \ee^{V(b_{N-1})-V(m_N)+C}
\leq c_8 \ee^{-H^+_{N-1}}
\end{equation}
(compare with (\ref{geom})).

Recall (\ref{cst}) that
$E_\omega (Z_1) \leq \delta^{-1} \Lambda_N$  and further that
$\Lambda_N \geq 1$ (see (\ref{Lambda>c})). Hence,
with $\kappa \geq 2(\delta^{-1}+1)+1$, adapting the derivation of (\ref{M>ell})
we get for $k_r = 2^r$, any $\ell \geq 1$ and $i \in [b_N,\eta_N)$
the bound
\begin{eqnarray*}
&& P_\omega (\exists n \in [T(i),T(i+1)), \;
L(n,N) > \kappa [\oo{K} + \oo{M}] \Lambda_N) 
\nonumber \\    
&\le& P_\omega(K < \ell) + P_\omega (Z_0 \geq \ell \Lambda_N) +
\sum_{r=0}^\infty P_\omega \Big( 
\sum_{j=1}^{k_r} (Z_j - E_\omega(Z_j)) \ge [\ell + k_r] \Lambda_N \Big)
\\
&=:&
 I_1 + I_2 + I_6(i) \,.
\nonumber 
\end{eqnarray*}
This applies also for $i \in [m_N,b_N)$, upon setting $I_6(i)=0$. Fixing
$\varepsilon > 0$ we take care of the term $I_1$ by choosing 
$\ell:= \lceil p(b_{N-1},m_N)^{-1}
\exp(-\frac{1}{3}(H_{N-1}^+)^{1-\varepsilon}) \rceil$
(see (\ref{4-1})). By (\ref{eq:p-bd}) such choice also handles the term $I_2$
(compare with (\ref{4-2})). All that remains is to deal with 
the sum of $I_6(i)$ over $[b_N,\eta_N)$. To this end, 
adapting the derivation of (\ref{eq:I3}), we get the bound
\begin{equation}\label{eq:I3i<eta}
I_6(i) 
\leq \hbox{Var}_\omega (Z_1) \sum_{r=0}^\infty \frac{k_r}{(\ell+k_r)^2}
\leq \frac{c_{14} \hbox{Var}_\omega (Z_1)}{\ell} \,.
\end{equation}
Recall the bound (\ref{eq:bd-var}) on $\hbox{Var}_\omega(Z_1(i))$
for $i \geq b$, the monotonicity of 
$i \mapsto U(b,i)$ and the fact that $P$-almost surely
$U(b,\eta_N) \leq H_{N-1}^+ - (H_{N-1}^+)^{1-\varepsilon}$.
Together with (\ref{eq:I3i<eta}), our choice of $\ell$ and 
the bound (\ref{eq:p-bd}), this results with 
$$
\sum_{i=b_N}^{\eta_N-1} I_6(i) \leq c_{15} m_{N+1}^4 
\ee^{-\frac{2}{3}(H_{N-1}^+)^{1-\varepsilon}} \,,
$$
holding for all $N$ large enough. As usual, 
by (\ref{bk<Hk2}) and Lemma \ref{l:Hk}, 
this concludes the proof.
\hfill$\Box$

\bigskip

\noindent {\it Proof of Proposition \ref{l:ratio}.} 
Our claim (\ref{ratio}) amounts to 
having $\p$-almost surely for $N$ large,
$$
\max_{n \in [T(m_N),\; T(m_{N+1}))} \, {L(n, N-1)\over \xi(n,
b_{N-1})} \le c_3 \, \Lambda_{N-1} \; ,
$$
which in view of Lemma \ref{l:Hk}
follows from (\ref{L>xi(b)}) by the Borel--Cantelli lemma.
Similarly, since $n \mapsto \xi(n,x)$ is monotone,
combining Lemma \ref{l:xi+xi} and (\ref{L>xi(b)}) we find 
by Lemma \ref{l:Hk} that 
$$
\sum_N P_\omega \Big(\exists n \geq T(m_N),\, 
 L(n,N) > c_4 \, [\,
\xi(n, b_{N-1}) + \xi(n, b_N) \, ] \Lambda_N \Big) < \infty ,
$$
for $c_4 = \max(c_3,\kappa)$. Applying the Borel--Cantelli lemma, this 
obviously implies (\ref{ratioN}).\hfill$\Box$

\subsection{
The effective width of the valleys}
\label{subs:expo}

\medskip
We consider next the asymptotic growth of the effective 
width $\Lambda_k$ (see (\ref{Lambda})), of the valleys.
\begin{proposition}
 \label{l:expo} 
 There exist constants $0< \gamma_- \le \gamma_+ <\infty$ such that
 \begin{equation}
    \gamma_- \le \limsup_{N\to \infty} \, {1\over \log N} \max_{1\le
    k\le N} \Lambda_k \le \gamma_+\, , \qquad \hbox{\rm $P$-a.s.}
    \label{expo}
 \end{equation}
\end{proposition}

\medskip

\noindent {\it Proof.}
We start by proving the lower bound in (\ref{expo}).
To this end, consider the events
$$
E_k := \left\{ V(b_k-i)-V(b_k)
\le c_{16}, \; \; \forall \, 0\le i\le c_{17} \, \log k\right\},
$$
for finite positive constants $c_{16}$ and $c_{17}$ to be chosen later.
Recall that $b_k - m_k \ge H_k^-/\log ({1-\delta\over \delta})$ for
the constant $\delta$ of (\ref{omega0}) and
$P$-almost surely $\log H_k^- \sim k$ for all large $k$
(by Lemma \ref{l:Hk}). Consequently, the interval
$[b_k - c_{17} \log k,b_k]$ lies inside the $k$-th valley
for all $k$ large enough, in which case the event $E_k$
implies that $\Lambda_k \ge c_{17} \ee^{-c_{16}} \log k$.
Since $E_k$ is adapted to the filtration
${\cal G}_k := \sigma\{ V(i), \; 0\le i\le \theta_{k+1} \}$, if
\begin{equation}
    \sum_k P \left( E_k \, | \, {\cal G}_{k-1} 
    \right) = \infty, \qquad \hbox{\rm $P$-almost 
    surely,}
    \label{levy-bc}
\end{equation}
then by L\'evy's Borel--Cantelli lemma (see, for example,
\cite[Page 518]{shiryaev}), we have that $P$-almost surely
$E_k$ occurs for infinitely many $k$, and therefore
$$
\limsup_{k\to \infty} \, {\Lambda_k \over \log k} \ge
c_{17}\ee^{-c_{16}} , \qquad \hbox{\rm $P$-a.s.}.
$$
This clearly yields the lower bound in (\ref{expo}), with
$\gamma_- = c_{17}\ee^{-c_{16}} >0$.

Turning to prove (\ref{levy-bc}), define, for any $\rho>0$,
\begin{eqnarray}
    \eta(\rho)
 &:=& \inf \Big\{ i > 0: \; V_k (i) - \min_{0\le j\le i} V_k (j) \ge
    \rho \Big\} ,
    \label{eta(n)}
    \\
    b(\rho)
 &:=& \sup \Big\{ i < \eta(\rho) : \; V_k (i) 
= \min_{0\le j\le \eta(\rho)} V_k (j)
    \Big\},
    \label{b(n)}
\end{eqnarray}
and the associated events
$$
E(\rho,k):= \{ V_k (b(\rho) - i)-V_k (b(\rho)) \le c_{16},
\; \; \forall \, 0\le i\le c_{17} \, \log k \}\,,
$$
where $(V_k(i):=V(i+\theta_k)-V(\theta_k), i \in \z_+)$ 
has the same law as $(V(i),i \in \z_+)$.
Recall that $H_{k-1}^+$, $\theta_k$ and $V(\theta_k)$ are 
${\mathcal G}_{k-1}$-measurable while 
$\eta_k-\theta_k =\eta(\rho)$ and $b_k-\theta_k = b(\rho)$ 
for $\rho=H_{k-1}^+$. Thus, $P(E_k|{\mathcal G}_{k-1})=P(E(\rho,k))$
for this choice of $\rho$. Since  $P$-almost surely,
$\log H_{k-1}^+ \sim k$ for $k\to \infty$ (by Lemma \ref{l:Hk}),  
the proof of (\ref{levy-bc}) is reduced to showing that 
for some $c_{18}>0$ and all $k$ large enough, 
\begin{equation}
    \inf_{\rho\ge \ee^{k/2}} P( E(\rho,k) ) \geq \frac{c_{18}}{k} \,.
    \label{levy-bc1}
\end{equation}
To verify (\ref{levy-bc1}), recall that by 
assumption (\ref{omega0}) the increments of the 
random walk $V$ are within $[-C,C]$ for some
$C=C(\delta)$ finite and positive. Further, (\ref{nondet}) yields
that $p_*:= \min \{ P(V(1) > 0) , \, P(V(1) \leq -2C/d)\} >0$ 
for some finite positive integer $d \geq 2$.
Restricting the  
increments of the walk $V(i+1)-V(i)$, $i \leq j-1$ 
to be strictly positive if $V(i) \leq C$ and at most $-2C/d$ otherwise, 
followed by $d$ increments which are at most $-2C/d$ each,
results in a sample of length $j+d$ for which the event  
$$
F_j := \bigcup_{\ell=j+1}^{j+d} \big\{ \, V(0)=0, \; V(i) \in (0,2C], 
\; \, 0\le i\le \ell-1, \; V(\ell) \in (-C,0] \, \big\}  
$$
holds. Hence, $P(F_j) \geq p_*^{j+d}$. In particular, if $c_{17}>0$
is small enough, then $P(F_j) \geq \frac{c_{18}}{k}$, 
for $j = \lceil c_{17} \log k \rceil$,  some $c_{18}>0$
and all $k$. Setting $V(\cdot)$ for $V_k(\cdot)$,
it is well known that the path $(V(i),\,0 \leq i \leq \eta(\rho))$
can be constructed by concatenating i.i.d. excursions of the 
walk $V(\cdot)$, each starting at $0$ and terminating at 
the first exit time of $(0,\rho)$. Then, 
$\eta(\rho)$ is the terminal time of the first excursion to exit 
via $[\rho,\infty)$ with $b(\rho)$ its starting time, and in concatenating the 
excursions en-route to the path, one adds to each excursion the 
(non-positive) values of all terminal points of preceding excursions. 
Adopting this construction, the event $E(\rho,k)$ 
occurs for $c_{16}=3C < \rho$ if the
last of the excursions which exit via $(-\infty,0]$ is 
in $F_j$ for $j = \lceil c_{17} \log k \rceil$, yielding the
bound (\ref{levy-bc1}) in view of the independence of these excursions.  
 
Turning to show the upper bound in (\ref{expo}),
let $\eta=\eta(\rho)$ and $b=b(\rho)$ be as in (\ref{eta(n)})--(\ref{b(n)}).
In the sequel we show that for some positive finite constants 
$c_{19}$, $c_{20}$ and $r_0 \geq 1$,
\begin{equation}
    \sup_{\rho \ge K_0} P\Big( \sum_{i=0}^{\eta-1} \ee^{- [V(i)-V(b)]} >
    c_{20} r \Big) \le \ee^{- c_{19} r}, \qquad \forall \, r \geq r_0.
    \label{expo1}
\end{equation}
Recall that conditional upon ${\mathcal G}_{k-1}$ 
the joint law of $V(i)-V(b_k)$ for $i \in [\theta_k,\eta_k]$ is the same
as the unconditional joint law of $V(i)-V(b(\rho))$ 
for $i \in [0,\eta(\rho)]$ upon taking 
$\rho=H_{k-1}^+$ (which is measurable on ${\mathcal G}_{k-1}$). Since
$H_{k-1}^+ \geq H_0^+ \geq K_0$, it thus follows from (\ref{expo1}) that    
$P( \sum_{i=\theta_k}^{\eta_k-1} \ee^{- [V(i)-V(b_k)]} > c_{20} r )
\le \ee^{- c_{19}\, r}$, for any $k \geq 1$. So, by the
Borel--Cantelli lemma, for some $\gamma_+<\infty$ and  
$P$-almost surely for $N$ large enough,
$$
\sum_{i=\theta_N}^{\eta_N-1} \ee^{-[\, V(i)-V(b_N)\,]} \le \gamma_+ \, \log N.
$$
\noindent Further, from the definition of $\theta_N$ and $H_{N-1}^+$ 
we know that $P$-almost surely for $N$ large enough, 
$$
\sum_{i=m_N}^{\theta_N-1} \ee^{-[\, V(i)-V(b_N)\,]} \le
\sum_{i=m_N}^{\theta_N-1} \ee^{- H_{N-1}^+} \le b_N
\, \ee^{- H_{N-1}^+} \le \ee^{-N} 
$$
(the last inequality being a consequence of (\ref{bk<Hk2})
and Lemma \ref{l:Hk}). Also, by (\ref{depth3}), for any
$\varepsilon>0$ and $P$-almost surely for all large $N$, 
$$ 
\sum_{i=\eta_N}^{m_{N+1}-1} \ee^{-[\, V(i)-V(b_N)\,]} \le
m_{N+1} \ee^{-(H_{N-1}^+)^{1-\varepsilon}} \le \ee^{-N} .
$$
Thus, we have that $P$-almost surely 
$$
\Lambda_N = \sum_{i=m_N}^{m_{N+1}-1}\ee^{-[\, V(i)-V(b_N)\,]} \le
2 \ee^{-N} + \gamma_+ \log N,
$$
for all large $N$, clearly yielding the upper bound in (\ref{expo}).

To complete the proof of the proposition, it thus 
remains only to prove (\ref{expo1}). To this end, setting 
$\eta=\eta(\rho)$ and $b=b(\rho)$, we consider the 
random variables $L(j) := \# \{ i < \eta: \; V(i) -
V(b) \in [j, j+1)\}$, $j\in \z_+$ (which depend on $\rho$ 
via $\eta$ and $b$) and the events
$$
A_{j,m} := \left\{ -(m+1) < V(b) \le -m, \; L(j) > c_{21} \ee^{j/2} r \right\},
$$
for $j,m \in \z_+$ and $c_{21}<\infty$ to be determined in the sequel.
Since $\{ L(j) > c_{21} \ee^{j/2} r \}$ is the disjoint union of 
$A_{j,m}$ and
$\sum_{i=0}^{\eta-1} \ee^{- [V(i)-V(b)]}
\le \sum_{j=0}^{\infty} \ee^{-j} L(j)$, it follows that 
\begin{equation}
    P\Big( \sum_{i=0}^{\eta-1} \ee^{- [V(i)-V(b)]} > c_{20} r \Big)
    \le \sum_{j=0}^{\infty} P\big( L(j) > c_{21} \, \ee^{j/2} r \big)
    = \sum_{j=0}^{\infty} \sum_{m=0}^\infty P(A_{j,m}) \,,
    \label{Ajm}
\end{equation}
provided $c_{20} \geq c_{21} \sum_{j=0}^\infty \ee^{-j/2}$.  

We thus proceed to bound $P(A_{j,m})$ for all $j,m$ and $\rho \geq K_0$.
To this end, as $V(\cdot)$ is a non-degenerate random 
walk of zero mean and bounded increments, for large  
positive integer $c_{22}$ we have that 
\begin{equation}
    q_* := \sup_{j\ge 0} P\big(\inf_{i\in [0, \, (j+c_{22})^2)} V(i) >
    -(j+2) \big) < 1.
    \label{proba<1}
\end{equation}
Next, fixing $j \in \z_+$, let $g=g(j)=(j+c_{22})^2 \geq 1$ and 
$R = R(j) = \lceil c_{21} \ee^{j/2} r/g(j) \rceil -1$, where 
$c_{21}$ is taken sufficiently large so that $R(j) \geq 1$ for any 
$j \in \z_+$ and $r \geq 1$. Fixing also $m \in \z_+$,
we consider the stopping times
\begin{eqnarray*}
    T_0 
 &:=& \inf \left\{ i \geq 0: \; V(i) \in (j-m-1, j-m+1) \right\},
    \\
    T_\ell 
 &:=& \inf \left\{ i \geq T_{\ell-1} + g : \; V(i) \in (j-m-1,
    j-m+1) \right\}, \qquad \ell \ge 1,
\end{eqnarray*}
and the associated stopped $\sigma$-fields ${\mathcal F}_\ell$.
Suppose the event $A_{j,m}$ holds.
Then, the random walk $(V(i), \; i \leq \eta-1)$
hits the interval $(j-m-1, j-m+1)$ more than 
$\lfloor c_{21} \ee^{j/2} r \rfloor \ge R g$ times, 
and hence $T_R < \eta$. In particular, as the walk $V(\cdot)$
can not reach $[\rho,\infty)$ for $i < \eta(\rho)$ and the event
$\Gamma_0 := \{ T_0 < \eta \}$ must hold as well, it follows that
$A_{j,m}$ is an empty set whenever $j-m-1 \geq \rho$.
Further, if $A_{j,m}$ holds, then by the preceding discussion also the 
events
$$
\Gamma_{\ell} := \big\{ \inf_{i \in [0,g)} \, V(T_{\ell-1}+i) > -(m+1) \big\}
$$
hold for $\ell=1,\ldots,R$. Finally, if $A_{j,m}$ holds then
$V(\eta) \geq \rho + V(b) > \rho - (m+1)$, while $V(i)>-(m+1)$
for all $i \in (T_R,\eta]$, implying that the event
$$
\Gamma_* = \big\{ \, V(T_R+i), i \geq 0,  \hbox { exits } (-(m+1),\rho-(m+1)] 
 \hbox { upwards} \big\}\,,
$$
holds as well. To summarize, we have seen that
$$
A_{j,m} \subseteq \Gamma_0 \cap \bigcap_{\ell=1}^R \Gamma_\ell 
\cap \Gamma_* \,.
$$
The following bounds apply
\begin{eqnarray}\label{eq:gstar}
P(\Gamma_* \,|\, {\mathcal F}_R) &\leq& \sup_{x \in (j,j+2)}
P \Big( V(\cdot) \hbox{ exits } (0,\rho] \hbox { upwards }  
\; \big| \; V(0) = x \Big) \leq \frac{j+2+C}{\rho+C} \,, \\
\label{eq:gl}
P(\Gamma_\ell \;|\; {\mathcal F}_{\ell-1}) &\leq&
\sup_{x \in (j-m-1,j-m+1)} 
P\Big( \inf_{i\in [0, g)} V(i) > -(m+1) \, \big| \, V(0)=x \Big) \leq q_* \,,
\end{eqnarray}
using (\ref{proba<1}) in the latter bound. Further, if 
$j-m+1 \leq -J \rho$ for some $J \in \z_+$, then 
considering the first downward crossing of $-k \rho$ for $k=1,\ldots,J$, 
leads to 
$$
P(\Gamma_0) = P( T_0 < \eta) \le 
P\Big(\, V(\cdot) \hbox{ exits } (-\rho+C,\rho) \hbox { downwards } \Big)^J  
\le \Big(\frac{\rho+C}{2\rho}\Big)^J  \,.
$$
Since $\Gamma_0$ is empty for $j-m-1 \geq \rho$, this implies 
that for some finite $c_{23}$ and positive $c_{24}$, 
\begin{equation}\label{eq:g0}
P(\Gamma_0) \leq c_{23} \ee^{- c_{24} |m-j|/\rho} \,\,, \qquad
\qquad \forall \rho \geq K_0, \;  j,m \in \z_+ \,.  
\end{equation}
With $\Gamma_\ell$ measurable on ${\mathcal F}_\ell$, upon applying the
strong Markov property at the stopping times $T_\ell$, $\ell=0,\ldots,R$, 
we get from (\ref{eq:gstar}), (\ref{eq:gl}) and (\ref{eq:g0}) that
\begin{eqnarray*}
P\left( A_{j,m} \right) &\le&
P (\Gamma_0\cap \bigcap_{\ell=1}^R \Gamma_\ell \cap \Gamma_*)  
\leq
c_{23} \ee^{- c_{24}|m-j|/\rho} \,\, q_*^R \,\, \frac{j+2+C}{\rho+C} \\
&\le& c_{25} \, {j+1\over \rho} \exp \Big[ -
c_{26}\, \Big( {|m-j| \over \rho} + {\ee^{j/2} r \over (j+c_{22})^2}
\Big) \Big] \,.
\end{eqnarray*}
This implies for some finite $c_{27}$ and all $\rho \geq K_0$,
$$
\sum_{m=0}^{\infty} P (A_{j,m}) \le c_{27} \, (j+1)
\exp\big(- c_{26} {\ee^{j/2} r \over (j+c_{22})^2}\big) \,.
$$
Plugging the latter bound into (\ref{Ajm}) yields (\ref{expo1}), thus 
concluding the proof of the proposition.
\hfill$\Box$

\section{Proof of Theorem \ref{t:main}}
\label{s:proof}

We start with a preliminary result.

\medskip

\begin{lemma}
\label{l:favsite}
We have  
$$
\lim_{n\to \infty} \, {\xi(n, b_{N_n}) + 
\xi(n, b_{N_n-1})
\over \max\limits_{y< 
m_{N_n-1}
}\xi(n, y)} = \infty, \qquad
\hbox{\rm $\p$-a.s.}
$$

\medskip

\end{lemma}
\noindent {\it Proof.} According to Proposition \ref{l:ratio},
we have $\p$-a.s. for $n$ large enough 
$$
\xi(n, b_{N_n}) + 
\xi(n, b_{N_n-1})
\ge
\max\left\{ {L(n, N_n)\over c_4 \, \Lambda_{N_n}}, \; {L(n,
N_n-1)\over c_3 \, \Lambda_{N_n-1}}\right\} \, .
$$ 
Recall that by Proposition \ref{l:expo}, $\p$-a.s. for all $n$ large enough
$$
\max\{\Lambda_{N_n-1}, \Lambda_{N_n}\}
\le 2 \gamma_+ \log N_n \,, 
$$
and by 
(\ref{negli}), for any $\delta<1$, also $\p$-a.s. 
\begin{eqnarray*}
    L(n, N_n-1)+L(n, N_n) 
 &\ge& \exp\left((\log n)^\delta\right) \sum_{1\le k< N_n-1} L(n,k)
    \\
 &\ge& \exp\left((\log n)^\delta\right)
    \max\limits_{y< m_{N_n-1}}\xi(n, y)\, .
\end{eqnarray*}
Hence, $\p$-a.s. for $n$ large enough,
$$
\xi(n, b_{N_n}) + 
\xi(n, b_{N_n-1})
 \ge
c_{28}\, (\log N_n)^{-1} \exp\left((\log
n)^\delta\right)\max\limits_{y< m_{N_n-1}}\xi(n, y) .
$$
Since $N_n \sim \log\log n $ for $n\to \infty$ (see
(\ref{N})), this proves the claim of the lemma.\hfill$\Box$

\bigskip

The rest of the section is devoted to the proof of Theorem \ref{t:main}.

\bigskip

\noindent {\it Proof of Theorem \ref{t:main}.} According to a
0--1 law in \cite{mama}, there exists a possibly degenerate constant
$c \in [0,\infty]$ such that
$$
\liminf_{n\to \infty} \, {\xi^*(n) \over n/\log\log\log n}
=c, \qquad \hbox{\rm $\p$-a.s.}
$$

\noindent (Though the 0--1 law was proved in \cite{mama} for
transient random walk in random environment, its proof remains
valid for our recurrent walk, with a reflecting barrier at the  
origin.)

It remains to check that $0<c<\infty$.

We start by showing that $c$ is positive. From  
Proposition \ref{l:ratio} we have that $\p$-a.s. for $n$ large enough,
$$
L(n,\, N_n-1) + L(n,\, N_n) \le (c_3+c_4)
[\xi(n,\, b_{N_n-1}) + \xi(n,\, b_{N_n})] \max_{k \leq N_n} \Lambda_k
$$
Hence, combining (\ref{N}) with the upper bound in
Proposition \ref{l:expo}, we have
$$
\liminf_{n\to \infty} \, {[\, \xi(n,\, b_{N_n-1}) + \xi(n,\,
b_{N_n})\, ]\, \log\log\log n \over L(n,\, N_n-1) + L(n,\, N_n)} \ge
{1\over (c_3+c_4) \gamma_+}\, , \qquad \hbox{\rm $\p$-a.s.}
$$
Since $\xi(n,\, b_{N_n-1}) + \xi(n,\, b_{N_n}) \le
2\xi^*(n)$ and $\p$-almost surely 
$n^{-1} (L(n,\, N_n-1) + L(n,\, N_n)) \to 1$
for $n \to \infty$ (as a consequence of Theorem \ref{t:spend}), this implies
that
$$
\liminf_{n\to \infty} \, {\xi^*(n)\, \log\log\log n \over n} \ge
{1\over 2(c_3+c_4) \gamma_+}\, , \qquad \hbox{\rm $\p$-a.s.}
$$
Consequently, $c \geq 1/(2 (c_3+c_4) \gamma_+) > 0$ as claimed.

Turning to show that $c<\infty$, note that  
if $n=T(m_N)$ then $N_n=N$ while $\xi(n,b_N)=0$.
Thus, by Lemma \ref{l:favsite}, $\p$-a.s, if $n=T(m_N)$ for $N$ large enough,
then $\xi^*(n)=\max_{x \in [m_{N-1},m_N)} \xi(n,x)$. Consequently, 
by (\ref{ratio-lb}), and the trivial inequality $L(T(m_N), N-1) \le
T(m_N)$, we have that
$$
    \limsup_{N\to \infty} \, { \xi^*(T(m_N)) \Lambda_{N-1} \over T(m_N)} 
    \le {1\over c_5}\, , \qquad \hbox{\rm
    $\p$-a.s.}
$$
By the lower bound in Proposition \ref{l:expo}, it follows that 
$\limsup_k \, (\log k)^{-1} \Lambda_{k-1} \ge \gamma_-$. Consequently 
$$
\liminf_{N\to \infty} \, {\xi^*(T(m_N)) \log N \over T(m_N)} 
\le {1\over c_5\, \gamma_- }\, , \qquad
\hbox{\rm $\p$-a.s.}
$$
Since $P$-almost surely $\log\log m_N \sim \log N$ for all $N$ large
enough (see (\ref{loglogm_k})), and $\p$-almost surely 
$\log\log T(x) \sim {1\over 2}\, \log x$ for $x\to \infty$ 
(see Fact \ref{f:revesz}), it follows that $\p$-almost surely 
$\log\log\log T(m_N) \sim \log N$ for $N\to \infty$. Therefore, 
$$
\liminf_{n\to \infty} \, {\xi^*(n) \log\log\log n \over n} \le
{1\over c_5\, \gamma_-}\, , \qquad \hbox{\rm $\p$-a.s.}
$$
We deduce that $c \leq 1/(c_5 \gamma_-)$ is finite 
and hence conclude the proof of Theorem \ref{t:main}.\hfill$\Box$

\bigskip
\bigskip
\bigskip

\noindent {\Large\bf Acknowledgement}

\medskip

\noindent We would like to thank an anonymous referee for helpful comments on the first version of the paper.

\bigskip
\bigskip


{\footnotesize 

\baselineskip=12pt

\noindent 
\begin{tabular}{lll}
& \hskip30pt Amir Dembo
    & \hskip60pt Nina Gantert \\ 
& \hskip30pt Department of Mathematics 
    & \hskip60pt Fachbereich Mathematik und Informatik \\   
& \hskip30pt Stanford University
    & \hskip60pt Universit\"at M\"unster \\   
& \hskip30pt Stanford, CA 94305
    & \hskip60pt Einsteinstrasse 62 \\   
& \hskip30pt U.S.A.
    & \hskip60pt D-48149 M\"unster \\  
& \hskip30pt {\tt amir@math.stanford.edu}
    & \hskip60pt Germany \\ 
& \phantom{-} 
    & \hskip60pt {\tt gantert@math.uni-muenster.de}
    \\ \\ 
& \hskip30pt Yuval Peres
    & \hskip60pt Zhan Shi \\ 
& \hskip30pt Department of Statistics
    & \hskip60pt Laboratoire de Probabilit\'es et Mod\`eles
      Al\'eatoires \\    
& \hskip30pt University of California Berkeley
    & \hskip60pt Universit\'e Paris VI \\   
& \hskip30pt 367 Evans Hall
    & \hskip60pt 4 place Jussieu \\   
& \hskip30pt Berkeley, CA 94720-3860
    & \hskip60pt F-75252 Paris Cedex 05 \\  
& \hskip30pt U.S.A. 
    & \hskip60pt France \\ 
& \hskip30pt {\tt peres@stat.berkeley.edu} 
    & \hskip60pt {\tt zhan@proba.jussieu.fr}
\end{tabular}

}

\end{document}